\documentclass[11pt]{article}

\usepackage{authblk}

\usepackage[utf8]{inputenc}
\usepackage{etoolbox}
\usepackage{xparse}
\usepackage[utf8]{inputenc} 
\usepackage[T1]{fontenc} 
\usepackage{amssymb,amsfonts,amsopn,dsfont,amscd}
\usepackage{amsthm}
\usepackage{latexsym}
\usepackage{bm}
\usepackage{mathtools}
\usepackage{nicefrac}
\usepackage{mathrsfs}
\usepackage{empheq}
\usepackage{enumitem}
\usepackage{pifont}
\usepackage{color}
\usepackage{psfrag}
\usepackage{tikz}
\usetikzlibrary{arrows}
\usetikzlibrary{positioning}
\usepackage{fancybox}
\usepackage{tcolorbox}
\usepackage{listings}
\usepackage{lstautogobble}
\usepackage{a4}
\usepackage{ifthen}
\usepackage{helvet}     
\usepackage{times}      
\usepackage{listings}
\usepackage{subcaption}
\usepackage{algorithm}
\usepackage{algpseudocode}
\usepackage{algpascal}
\usepackage{esvect}
\usepackage{xcolor}

\usepackage[percent]{overpic}

\usepackage[final]{showkeys}
\usepackage[hidelinks]{hyperref}
\usepackage[T1]{fontenc}

\usepackage{musicography}
\usepackage[noabbrev, capitalize]{cleveref}

\newcommand{\iEG}{\mathcal{I}^{\operatorname{av}}_{h0}}

\newcommand{\R}{\mathbb{R}}

\newcommand{\X}{\mathbf{X}}
\newcommand{\Curl}{\nabla \times }
\newcommand{\Grad}{\nabla}
\newcommand{\Div}{\nabla \cdot}

\newcommand{\bgamma}{\boldsymbol{\gamma}}
\newcommand{\ttrace}{\boldsymbol{\gamma}_{\parallel}}
\newcommand{\rtrace}{\boldsymbol{\gamma}_t}
\newcommand{\ntrace}{\boldsymbol{\gamma}_{\Curl}}
\newcommand{\blambda}{\boldsymbol{\lambda}}

\newcommand{\VectorLtwo}{\mathbf{L}^2(\Omega)}
\newcommand{\VectorLtwoGamma}{\mathbf{L}^2(\Gamma)}
\newcommand{\VectorLtwoGammaPar}{\mathbf{L}^2_{\parallel}(\Gamma)}
\newcommand{\VectorH}{\mathbf{H}}
\newcommand{\Hcurl}{\VectorH(\operatorname{curl}, \Omega)}
\newcommand{\Hcurlz}{\VectorH_0(\operatorname{curl}, \Omega)}
\newcommand{\Hdiv}{\VectorH(\operatorname{div}, \Omega)}
\newcommand{\Hzdiv}{\VectorH_0(\operatorname{div}, \Omega)}

\newcommand{\VectorHs}{\VectorH^s(\Omega)}

\newcommand{\VectorHtwo}{\VectorH^2(\Omega)}

\newcommand{\HoneAvgZero}{H_\ast^1(\Omega)}

\newcommand{\V}{\mathbf{V}}
\newcommand{\Vh}{\mathbf{V}_h}
\newcommand{\Vhz}{\mathring{\mathbf{V}}_h}
\newcommand{\Qh}{Q_h}

\let\div\relax
\DeclareMathOperator{\div}{div}

\DeclareMathOperator{\curl}{curl}

\newcommand{\bfu}{\mathbf{u}}
\newcommand{\bfv}{\mathbf{v}}
\newcommand{\bfw}{\mathbf{w}}
\newcommand{\bfz}{\mathbf{z}}
\newcommand{\bff}{\mathbf{f}}
\newcommand{\bfg}{\mathbf{g}}

\newcommand{\nvec}{\mathbf{n}}

\newcommand{\bzetau}{\mathbf{w}}
\newcommand{\bzetap}{\varphi}

\newcommand{\Pu}{P_{\Vhz}\bzetau}
\newcommand{\Pp}{P_{\Qh}\bzetap}

\newcommand{\bfzero}{\mathbf{0}}
\newcommand{\symgrad}{\boldsymbol{\epsilon}}

\renewcommand{\epsilon}{\ensuremath\varepsilon}

\renewcommand{\phi}{\ensuremath{\varphi}}

\newtheorem{theorem}{Theorem}
\newtheorem{remark}[theorem]{Remark}
\newtheorem{corollary}[theorem]{Corollary}
\newtheorem{lemma}[theorem]{Lemma}

\numberwithin{theorem}{section}
\numberwithin{equation}{section}

\hypersetup{pageanchor=false}

\title{H(curl)-based approximation of the Stokes problem with weakly enforced no-slip boundary conditions}

\author[1]{Wietse M. Boon}
\author[2]{Wouter Tonnon\thanks{Corresponding author: \href{mailto:wouter.tonnon@sam.math.ethz.ch}{wouter.tonnon@sam.math.ethz.ch}}}
\author[3]{Enrico Zampa}

\affil[1]{NORCE Norwegian Research Centre, Bergen 5008, Norway}
\affil[2]{SAM, ETH Zürich, CH-8092 Zürich, Switzerland}
\affil[3]{Department of Mathematics, University of Vienna, Vienna, Austria}

\date{\today}

\begin{document}

\maketitle

\begin{abstract}
    In this work, we show how to impose no-slip boundary conditions for an $\Hcurl$-based formulation for incompressible Stokes flow, which is used in structure-preserving discretizations of Navier-Stokes and magnetohydrodynamics equations. At first glance, it seems straightforward to apply no-slip boundary conditions: the tangential part is an essential boundary condition on $\Hcurl$ and the normal component can be naturally enforced through integration-by-parts of the divergence term. However, we show that this can lead to an ill-posed discretization and propose a Nitsche-based finite element method instead. We analyze the discrete system, establishing stability and deriving a priori error estimates. Numerical experiments validate our analysis and demonstrate optimal convergence rates for the velocity field.
\end{abstract}

\section{Introduction}

The equations governing incompressible Stokes flow arise as simplifications of various, distinct physical systems. Besides describing Navier-Stokes flow at low Reynold's numbers, the same equations arise as the incompressible limit of linearly elastic solid materials \cite[Rem.~8.1.2]{boffi2013mixed}. Moreover, the well-developed theory of Stokes discretizations has led to efficient numerical schemes for Biot poroelasticity \cite{rodrigo2018new}.

We consider the Stokes problem with the fluid velocity in $\Hcurl$, which is primarily motivated by two applications: the Navier–Stokes (NS) equations and magnetohydrodynamics (MHD). In the context of the NS equations, choosing the velocity in 
$\Hcurl$ enables the application of the efficient semi-Lagrangian method developed in \cite{TonnonHiptmairSemiLagrangianNavierStokes}. For MHD systems, this choice facilitates the discrete preservation of cross-helicity, an important invariant of the continuous problem, as shown in \cite{HuLeeXu21,LaakmanHuFarrell23, ShipengRuijie}. However, working in $\Hcurl$ introduces analytical and numerical challenges, particularly in handling boundary conditions and in the formulation of the continuous problem. To better understand and address these difficulties, we study a simplified setting: the Stokes problem.

The Stokes problem with $\Hcurl$-conforming velocity was considered in \cite{boon2024h} in the context of slip boundary conditions, in which only the normal component of the velocity is enforced to be zero. This work presents the next chapter in which we consider the case of no-slip boundary conditions. 

In this article, we show that imposing no-slip as essential boundary conditions leads to an ill-posed problem, despite the fact that tangential traces are well-defined in $\Hcurl$. We therefore propose and analyze a Nitsche-type method instead. We derive a priori error estimates in carefully chosen discrete norms that include a contribution of the boundary terms. The expected orders of convergence are validated by numerical experiments. 

The remainder of this article is organized as follows. \Cref{sec:VVP} introduces the Stokes problem and \Cref{sec:vf} presents the functional setting. Our Nitsche-based method is proposed and analyzed in \Cref{sec:discrete_vf}, and \Cref{sec:numerics} contains the numerical experiments. Conclusions are presented in \Cref{sec:concluding_remarks}.

\subsection{The Stokes problem in rotation form with no-slip boundary condition}
\label{sec:VVP}

Let $\Omega\subset \R^d$, $d = 2, 3$ be a bounded Lipschitz domain and denote by $\Gamma$ its boundary. Formally, we seek a velocity field $\bfu: \Omega\to\R^d$ and pressure field $p: \Omega\to\R$ such that
\begin{subequations}\label{eq:continuousSymGradStokes}
\begin{align}
    - 2 \Div\symgrad(\bfu) + \Grad p &= \mathbf{f}, & \text{ in }&\Omega, \label{eq:continuousSymGradStokesa} \\
    \Div \bfu &= 0,  & \text{ in }&\Omega, \\
    \bfu &= \bfg, & \text{ on }&\Gamma,
    \end{align}
\end{subequations}
where $\bff:\Omega\mapsto\R^d$ is a given forcing term and $\bfg:\Gamma\mapsto\R^d$ is a given boundary term. We will use $\nvec:\Gamma\mapsto\R^d$ to denote the outward oriented, unit vector normal to $\Gamma$. Moreover, $\symgrad(\bfu)$ denotes the symmetric gradient of $\bfu$, that is,
\begin{equation*}
    \symgrad(\bfu) \coloneqq \frac12 \left(\nabla \bfu + (\nabla\bfu)^T \right).
\end{equation*}

For sufficiently smooth $\bfu$ with $\Div\bfu=0$, we can rewrite the second order term in \Cref{eq:continuousSymGradStokesa} as
\begin{equation} \label{eq: rewrite symgrad}
    \begin{aligned}
        -2 \Div\symgrad(\bfu) &= -\Delta\bfu - \Grad (\Div\bfu)\\
        &= \Curl\Curl\bfu - 2\Grad(\Div\bfu)\\
        &= \Curl\Curl\bfu.
    \end{aligned}
\end{equation}
We aim to construct a variational formulation of the following form: seek $\bfu\in \V$ and $p\in Q$ such that
\begin{subequations}\label{eq:WeakContinuousCurlFormulation}
\begin{align}
    (\Curl\bfu,\Curl\bfv) +(\nabla p,\bfv) &= (\mathbf{f},\bfv),\label{eq:cont1}\\
    (\bfu,\nabla q) &= 0 \label{eq:cont2}
\end{align}
\end{subequations}
for all $\bfv\in \V$ and $q\in Q$, where $(\cdot, \cdot)$ denotes the $L^2(\Omega)$ product. Note that \cref{eq:cont2} yields $\bfu\cdot\nvec=0$. Thus, it remains to enforce $\bfu\times\nvec=\bfzero$. Since $\bfu\times\nvec=\bfzero$ is a well-defined trace for any $\bfu\in\Hcurl$, one might consider adding the condition $\bfu\times\nvec=\bfzero$ to $\Hcurl$ as an essential condition. However, we will see in \cref{sec:naiveapproach} that this leads to an ill-posed problem. Instead, we propose a Nitsche-type approach in \Cref{sec:discrete_vf}.

\section{The functional and topological framework}
\label{sec:vf}
In this section, we first describe the functional setting and introduce the notational conventions used throughout this work. 
\subsection{Functional Setting}
Let $\Omega\subset \R^d$, $d = 2, 3$, be a domain that is either a Lipschitz polygon or has a $C^{1,1}$-boundary and has a finite first Betti number. As mentioned in the introduction, we denote the $L^2(\Omega)$ inner product by parentheses and the analogous product on its boundary $\Gamma$ by angled brackets:
\begin{align*}
    (\bfu, \bfv) &\coloneqq \int_\Omega \bfu \cdot \bfv\, \mathrm{dx}, &
    \langle \bfu, \bfv \rangle &\coloneqq \int_\Gamma \bfu \cdot \bfv\,\mathrm{dS}.
\end{align*}
With a slight abuse of notation, we denote inner products for scalar functions using the same parentheses and brackets. We need the classical Sobolev spaces $W^{s, q}(\Omega)$ equipped with the Sobolev-Slobodeckij norm $\lVert \cdot \rVert_{s, \Omega, q} $, see e.g. \cite[Chapter 2.2.2]{EG1}. Then we set $H^s(\Omega) \coloneqq W^{s,2}(\Omega)$, and we denote its norm by $\lVert \cdot \rVert_{s, \Omega} \coloneqq \lVert \cdot \rVert_{s, \Omega, 2}$. Following this convention, $\lVert \cdot \rVert_{\Omega, q} \coloneqq \lVert \cdot \rVert_{0, \Omega, q}$ denotes the $L^q$ norm and $\lVert \cdot \rVert_{\Omega} \coloneqq \lVert \cdot \rVert_{\Omega, 2}$ denotes the $L^2$ norm. 
All these spaces can be defined also on $\Gamma$. We recall also the definitions of the following classical Hilbert spaces
\begin{align*}
    \HoneAvgZero &\coloneqq \{q\in H^1(\Omega)\mid (q,1) = 0\},\\
    \Hcurl &\coloneqq \{ \bfv \in \VectorLtwo \mid \Curl \bfv \in (L^2(\Omega))^{2d-3}\}, \\
     \Hdiv &\coloneqq \{ \bfv \in \VectorLtwo \mid \Div \bfv \in L^2(\Omega)\}, 
\end{align*}
and their associated norms
\begin{align} \label{eq: norms}
    \lvert q \rvert_{1, \Omega} &\coloneqq \lVert \nabla q \rVert_{\Omega}, \\
    \lVert\bfv\rVert_{\curl, \Omega}^2 &\coloneqq \lVert\bfv\rVert_\Omega^2 + \lVert\nabla\times\bfv\rVert_\Omega^2. \\
    \lVert\bfv\rVert_{\div, \Omega}^2 &\coloneqq \lVert\bfv\rVert_{\Omega}^2 + \lVert\nabla\cdot\bfv\rVert_{\Omega}^2.
\end{align}
We continue by defining $\Hcurlz$ and $\Hzdiv$ as the closure of $\boldsymbol{\mathcal{C}}_c^{\infty}(\Omega)$ with respect to the $\lVert \cdot \rVert_{\mathrm{curl}, \Omega}$ and $\lVert \cdot \rVert_{\mathrm{div}, \Omega}$ norms, respectively. In general, given a scalar-valued space $V$, we denote in bold $\mathbf{V}$ its vector-valued counterpart. In this work $a \lesssim b$ means $a \leq Cb$ where $C$ is a constant independent of the mesh size $h$ (when considering mesh-dependent quantities) and of $q$ (when considering $L^q$ spaces).

\subsection{Tangential traces}
We briefly recall some useful results from \cite{BuffaCostabelSheen}. We define
\begin{align*}
    \VectorLtwoGammaPar&\coloneqq \{ \boldsymbol{\zeta}\in \VectorLtwoGamma\mid \boldsymbol{\zeta}\cdot \nvec = 0\}
\end{align*}
Let $\bgamma: \mathbf{H}^1(\Omega) \to \mathbf{H}^{\frac{1}{2}}(\Gamma)$ be the standard trace operator and let $\bgamma^{-1}: \mathbf{H}^{\frac{1}{2}}(\Gamma)\to \mathbf{H}^1(\Omega)$ be one of its continuous right inverses, e.g.~the harmonic extension. 
We are now ready to define the tangential trace $\ttrace:\boldsymbol{\mathcal{C}}^{0}(\overline{\Omega})\to\VectorLtwoGammaPar$ and the tangential trace $\rtrace:\boldsymbol{\mathcal{C}}^{0}(\overline{\Omega})\to\VectorLtwoGammaPar$ operators as
\begin{align*}
    \ttrace(\bfv) &\coloneqq \nvec\times (\bfv |_{\Gamma}\times \nvec), & \rtrace(\bfv) & \coloneqq \bfv|_{\Gamma}\times \nvec, &
    \forall \bfv &\in \boldsymbol{\mathcal{C}}^{0}(\overline{\Omega}).
\end{align*}
In this work we will use only $\ttrace$, but we need $\rtrace$ to define the appropriate functional spaces on the boundary.
Using \cite[Thm.~3.10]{EG1}, $\ttrace$ and $\rtrace$ can be extended to bounded maps $\ttrace:\mathbf{H}^{1}(\Omega)\to \VectorLtwoGammaPar$ and $\rtrace: \mathbf{H}^{1}(\Omega)\to\VectorLtwoGammaPar$. With these operators we define 
\begin{equation*}
	\mathbf{H}^{\frac{1}{2}}_{\parallel}(\Gamma) \coloneqq \ttrace\circ \bgamma^{-1}(\mathbf{H}^{\frac{1}{2}}(\Gamma)), \qquad \mathbf{H}^{\frac{1}{2}}_{t}(\Gamma) \coloneqq \rtrace\circ \bgamma^{-1}(\mathbf{H}^{\frac{1}{2}}(\Gamma)).
\end{equation*}
with norms
\begin{align*}
	\lVert \blambda \rVert_{\frac{1}{2}, \parallel, \Gamma} & \coloneqq \inf_{\substack{
        \boldsymbol{\varphi} \in \mathbf{H}^{\frac{1}{2}}(\Gamma) \\
        \ttrace\circ \bgamma^{-1}(\boldsymbol{\varphi}) = \blambda}} 
        \lVert \boldsymbol{\varphi} \rVert_{\frac{1}{2}, \Gamma}, &
	\lVert \blambda \rVert_{\frac{1}{2}, t, \Gamma} & \coloneqq \inf_{\substack{
        \boldsymbol{\varphi} \in \mathbf{H}^{\frac{1}{2}}(\Gamma) \\
        \rtrace\circ \bgamma^{-1}(\boldsymbol{\varphi}) = \blambda}} 
        \lVert \boldsymbol{\varphi} \rVert_{\frac{1}{2}, \Gamma},
\end{align*}
Define $\mathbf{H}^{-\frac{1}{2}}_{\parallel}(\Gamma)$ and $\mathbf{H}^{-\frac{1}{2}}_t(\Gamma)$ as the dual spaces of $\mathbf{H}^{\frac{1}{2}}_{\parallel}(\Gamma)$ and $\mathbf{H}^{\frac{1}{2}}_t(\Gamma)$, respectively. The corresponding dual norms are respectively denoted by $\lVert \cdot \rVert_{-\frac{1}{2}, \parallel, \Gamma}$ and $\lVert \cdot \rVert_{-\frac{1}{2}, t, \Gamma}$.
\begin{lemma} 
\label{lemma:boundHmh}
    If $\boldsymbol{\varphi} \in \VectorLtwoGammaPar$, then $\lVert \boldsymbol{\varphi} \rVert_{-\frac{1}{2}, t, \Gamma} \lesssim \lVert \boldsymbol{\varphi}\rVert_{\Gamma}$. 
\end{lemma}
\begin{proof}
If $\blambda\in \mathbf{H}^{\frac{1}{2}}_{t}(\Gamma)$, then there exists $\widetilde{\blambda}\in\mathbf{H}^{\frac{1}{2}}(\Gamma)$ such that $\blambda = \rtrace \circ \bgamma^{-1}(\widetilde{\blambda})$. This implies that
\begin{equation*}
 	\lVert \blambda \rVert_{\Gamma} = \lVert \rtrace \circ \bgamma^{-1} (\widetilde{\blambda}) \rVert_{\Gamma} 
 	\lesssim \lVert \bgamma^{-1}(\widetilde{\blambda}) \rVert_{1 ,\Omega} 
 	\lesssim \lVert \widetilde{ \boldsymbol{ \lambda} } \rVert_{\frac{1}{2}, \Gamma}.
\end{equation*}
Since this is valid for any $\widetilde{\blambda}$ satisfying $\blambda = \rtrace \circ \bgamma^{-1}(\widetilde{\blambda})$, we obtain 
\begin{equation*}
	\lVert \blambda \rVert_{\Gamma} \lesssim \inf_{\substack{\widetilde{\blambda} \in\mathbf{H}^{\frac{1}{2}}(\Gamma) \\ \rtrace\circ\bgamma^{-1}(\widetilde{\blambda}) = \blambda }} \lVert \widetilde{\blambda} \rVert_{\frac{1}{2}, \Gamma} = \lVert \blambda \rVert_{\frac{1}{2}, t, \Gamma}.
\end{equation*}
As a consequence, if $\boldsymbol{\varphi} \in \VectorLtwoGammaPar$, it holds that
\begin{equation*}
	\begin{split}
		\lVert \boldsymbol{\varphi} \rVert_{-\frac{1}{2}, t, \Gamma} &= \sup_{\blambda\in\mathbf{H}^{\frac{1}{2}}_{t}(\Gamma)}\frac{\langle \boldsymbol{\varphi}, \blambda\rangle_*}{\lVert\blambda \rVert_{\frac{1}{2},t,\Gamma}}  \\ 
		&= \sup_{\blambda\in\mathbf{H}^{\frac{1}{2}}_{t}(\Gamma)}\frac{\langle \boldsymbol{\varphi}, \blambda\rangle_{\Gamma}}{\lVert\blambda \rVert_{\frac{1}{2},t,\Gamma}}  \\ 
		&= \sup_{\blambda\in\mathbf{H}^{\frac{1}{2}}_{t}(\Gamma)}\frac{\lVert \boldsymbol{\varphi} \rVert_{\Gamma} \lVert \blambda \rVert_{\Gamma} }{\lVert\blambda \rVert_{\frac{1}{2},t,\Gamma}} \\ 
		& \lesssim \lVert \boldsymbol{\varphi} \rVert_{\Gamma}.
	\end{split}
\end{equation*}
\end{proof}
On the other hand, it was shown in \cite{BuffaCostabelSheen} that $\ttrace$ can be extended to a bounded map $\ttrace:\Hcurl\to \mathbf{H}^{-\frac{1}{2}}_t(\Gamma)$.  We summarize the properties of $\ttrace$ in the following lemma.
\begin{lemma} \label{thm:TraceBounded}
    $\ttrace:\mathbf{H}^1(\Omega)\to \mathbf{H}^{\frac{1}{2}}_{\parallel}(\Gamma)$ and $\ttrace:\Hcurl\to \mathbf{H}^{-\frac{1}{2}}_t(\Gamma)$ are bounded.
\end{lemma}

\subsection{On topologically non-trivial domains}
We conclude this section by discussing some topological properties of the domains we consider. We denote by $\mathfrak{H}^1(\Omega)$ the associated space of Neumann harmonic vector fields:
\begin{equation*}
\begin{split}
    \mathfrak{H}^1(\Omega)  &\coloneqq \mathbf{H}(\operatorname{curl}_0, \Omega)\cap \mathbf{H}_0(\operatorname{div}_0, \Omega)\\ & \coloneqq \{ \mathbf{v}\in\mathbf{H}(\operatorname{curl}, \Omega)\mid \Curl \mathbf{v} = \bfzero\} \cap \{ \mathbf{v} \in \mathbf{H}_0(\operatorname{div}, \Omega) \mid \Div \mathbf{v} = 0\}.
    \end{split}
\end{equation*}
Its dimension equals the first Betti number, which we assume to be finite.

\begin{theorem} \label{thm:HarmonicFunctionsBound}
    Let $\Omega\subset \R^d$, $d=2,3$, be a Lipschitz domain. Then, there exists a $C>0$ such that for every $\mathfrak{h}\in \mathfrak{H}^1(\Omega)$
    \begin{align} \label{eq:L2BoundaryBoundContinuous}
        \|\mathfrak{h}\|_{\curl, \Omega} \leq C \|\ttrace(\mathfrak{h})\|_{-\frac{1}{2},t,\Gamma}.
    \end{align}
\end{theorem}
\begin{proof}
    We claim that $\lVert \ttrace(\cdot)\rVert_{\mathbf{H}^{-\frac{1}{2}}(\Gamma)}$ is a norm on $\mathfrak{H}^1(\Omega)$. In fact, let $\mathfrak{h}\in\mathfrak{H}^1(\Omega)$ satisfy $\lVert \ttrace(\mathfrak{h}) \rVert_{\mathbf{H}^{-\frac{1}{2}}(\Gamma)} = 0$. Then $\mathfrak{h}$ belongs to $\mathbf{H}_0(\operatorname{curl}_0, \Omega)\cap \mathbf{H}_0(\operatorname{div}_0, \Omega)$. In particular, $\mathfrak{h}\in\mathbf{H}^1_0(\Omega)$ by \cite[Theorem 2.5]{Amr+98}. Then \cite[Remark 2.6]{Amr+98} implies that $\mathfrak{h} = 0$, proving the claim. Then \eqref{eq:L2BoundaryBoundContinuous} follows from the equivalence of norms in finite dimensional spaces.
\end{proof}

\section{A Nitsche-based discretization method}
\label{sec:discrete_vf}

In this section, we show that a naïve approach to discretize
\cref{eq: rewrite symgrad} with the appropriate boundary conditions is not well-posed on all meshes. To avoid restriction on the mesh type and the relative technicalities in the analysis, we pursue a different strategy, proposing a Nitsche-type approach. 

\subsection{Preliminary definitions and inequalities}

Let $\{\Omega_h\}_h$ be a quasi-uniform and uniformly shape regular family of meshes on $\Omega$ with $h > 0$ the mesh size. Let $\Gamma_h$ denote the $(d-1)$-dimensional boundary mesh obtained by restricting $\Omega_h$ to $\Gamma$. 

Let $\{\Vh\}_h$ and $\{\Qh\}_h$ be asymptotically dense families of $\Omega_h$-piecewise-polynomial subspaces of $\Hcurl$ and $\HoneAvgZero$ respectively, satisfying
\begin{align*}
    \Grad \Qh \subseteq \Vh,
\end{align*} 
and the following approximation estimates:
\begin{align*}
    \inf_{\bfv_h \in \Vh}\lVert \bfu-\bfv_h \rVert_{\mathrm{curl}, \Omega}& \lesssim h^r (\lvert \bfu \rvert_{r,\Omega}+ \lvert \Curl \bfu\rvert_{r,\Omega}), \\
    \inf_{q_h\in Q_h} \lvert p - q_h \rvert_{1,\Omega} &\lesssim h^r\lvert p\rvert_{r+1, \Omega}
\end{align*}
for some positive integer $r$.
Additionally, we define $\X_h$ as the orthogonal complement of $\Grad\Qh$ in $\Vh$:
\begin{align*}
    \X_h \coloneqq \{ \bfv_h \in \Vh \mid (\bfv_h, \nabla q_h) = 0, \forall q_h \in Q_h \}.
\end{align*}
An important subspace of $\X_h$ is the space of discrete Neumann harmonic vector fields:
\begin{equation*}
    \mathfrak{H}^1_h \coloneqq \{ \mathbf{v}_h\in \X_h\mid \Curl \mathbf{v}_h = \bfzero\}.
\end{equation*}
Defining $\mathbf{Z}_h$ as the $L^2$-orthogonal complement of $\mathfrak{h}^1_h$ in $\mathbf{X}_h$, it is possible to show the following Hodge-Helmholtz orthogonal decomposition (see \cite[Thm.~4.5]{Arnold2018FiniteCalculus}):
\begin{equation}
    \mathbf{V}_h = \Grad Q_h \oplus^{\perp} \X_h = \Grad Q_h \oplus^{\perp} \mathbf{Z}_h \oplus^{\perp} \mathfrak{H}^1_h.
    \label{eq:HH_decomp}
\end{equation}
We define  $\mathring{Q}_h$ and $\Vhz$ be the subspaces of $\Qh$ and $\Vh$ with essential boundary conditions, respectively:
\begin{align*}
    \mathring{Q}_h &\coloneqq \Qh \cap H_0^1(\Omega), &
    \Vhz &\coloneqq \Vh \cap \VectorH_0(\curl, \Omega).
\end{align*}
The following trace operator $\ntrace:\Vh\mapsto \VectorLtwoGamma$ will also be needed
\begin{align*}
    \ntrace(\bfv_h) &\coloneqq \nvec\times (\Curl \bfv_h)|_\Gamma, & 
    \forall \bfv_h &\in \Vh.
\end{align*}

\begin{lemma}[Discrete trace inequalities] \label{lem: discrete trace ineq}
    There exist constants $C_n > 0$ and $C_\parallel > 0$ such that
    \begin{subequations}
    \begin{align}
        \| \ntrace(\bfv_h) \|_\Gamma &\leq C_n h^{-\frac12} \| \nabla \times \bfv_h \|_\Omega, \label{eq:trace_ineq_gamman}\\
        \| \ttrace(\bfv_h) \|_\Gamma &\leq C_\parallel h^{-\frac12} \| \bfv_h \|_\Omega, \label{eq:trace_ineq_gammapar}
            \end{align}
        \end{subequations}
    for each $\bfv_h \in \Vh$.
\end{lemma}
\begin{proof}
    Both inequalities follow from \cite[Lem.~12.8]{EG1}.
\end{proof}

\begin{lemma} \label{lem:L2BoundaryBoundDiscrete}
    There exists a $C>0$ such that for all $h>0$ small enough and every discrete Harmonic 1-form $\mathfrak{h}_h \in \mathfrak{H}^1_h(\Omega)$
    \begin{align*} 
        \|\mathfrak{h}_h\|_{\curl, \Omega} \leq C \|\ttrace(\mathfrak{h}_h)\|_{-\frac{1}{2},t,\Gamma}
    \end{align*}
    with $C>0$ independent of $h$.
\end{lemma}
\begin{proof}
    Take $\mathfrak{h}_h \in \mathfrak{H}^1_h(\Omega)$ arbitrary. Using \cite[Lem.~5.9]{AFW06} and \cite[Thm.~5.6]{AFW06}, there exists $\mathfrak{h}\in \mathfrak{H}^1(\Omega)$ such that $\lVert \mathfrak{h} \rVert_{\Omega} \leq \lVert \mathfrak{h}_h \rVert_{\Omega}$ and 
    \begin{align*}
        \lVert \mathfrak{h}-\mathfrak{h}_h\rVert_{\Omega} &\leq Ch^k \|\mathfrak{h}\rVert_{k, \Omega}, 
        & k &= \min\{ r, s \},
    \end{align*}
    in which $s>0$ is such that $\mathfrak{h} \in \VectorHs$ and $r \ge 1$ is the approximation order of the space $\Vh$ in $\Hcurl$. Note that since all norms on $\mathfrak{H}^1(\Omega)$ are equivalent, it also holds 
    \begin{align*}
        \lVert \mathfrak{h}-\mathfrak{h}_h\rVert_{\Omega} &\leq Ch^k \|\mathfrak{h}\rVert_{ \Omega}, 
        & k &= \min\{ r, s \},
    \end{align*}
    Moreover, since we assumed that $\Omega$ is a Lipschitz domain, \cite[Theorem 11.2]{MitreaMitreaTaylor} implies that $k\geq\frac{1}{2}$. 
    Using \cref{thm:HarmonicFunctionsBound} and \Cref{thm:TraceBounded}, we obtain
    \begin{align*}
        \|\mathfrak{h}_h\|_\Omega &\leq \|\mathfrak{h}\|_\Omega+\|\mathfrak{h}_h-\mathfrak{h}\|_\Omega\\
        &\leq C\|\ttrace(\mathfrak{h})\|_{-\frac{1}{2},t,\Gamma} + Ch^k \|\mathfrak{h}\rVert_{ \Omega}\\
        &\leq C\|\ttrace(\mathfrak{h}_h)\|_{-\frac{1}{2},t,\Gamma} +C\|\ttrace(\mathfrak{h}-\mathfrak{h}_h)\|_{-\frac{1}{2},t,\Gamma} +Ch^k \|\mathfrak{h}\rVert_{ \Omega}\\
        &\leq C\|\ttrace(\mathfrak{h}_h)\|_{-\frac{1}{2},t,\Gamma} +C\|\mathfrak{h}-\mathfrak{h}_h\|_{\curl, \Omega} + Ch^k \|\mathfrak{h}\rVert_{ \Omega}\\
        &= C\|\ttrace(\mathfrak{h}_h)\|_{-\frac{1}{2},t,\Gamma} +C\|\mathfrak{h}-\mathfrak{h}_h\|_\Omega + Ch^k \|\mathfrak{h}\rVert_{ \Omega}\\
        &\leq C\|\ttrace(\mathfrak{h}_h)\|_{-\frac{1}{2},t,\Gamma} +C^*h^k \|\mathfrak{h}\rVert_{\Omega}\\
        & \leq C\|\ttrace(\mathfrak{h}_h)\|_{-\frac{1}{2},t,\Gamma} +C^*h^k \|\mathfrak{h}_h\rVert_{ \Omega}.
    \end{align*}
    Choosing $h>0$ small enough such that $C^*h^k \leq \frac{1}{2}$ concludes the proof.
\end{proof}

\subsection{Failure of essential no-slip boundary conditions} 
\label{sec:naiveapproach}
Recall that $\Vhz$ and $\Qh$ are subspaces of $\Hcurlz$ and $\HoneAvgZero$, respectively. Then, we can define the following variational problem: seek $\bfu_h \in \Vhz$ and $p_h \in \Qh$ such that
\begin{subequations} \label{eq:wrongDiscreteFormulation}
    \begin{align}
        (\Curl \bfu_h,\Curl\bfv_h) + (\nabla p_h,\bfv_h) &= (\mathbf{f},\mathbf{v}_h),\\
        (\bfu_h,\nabla q_h)&= 0 \label{eq:wrongWeakDiv}
    \end{align}
\end{subequations}
for all $\bfv_h \in \Vhz$ and $q_h \in \Qh$. This means that we enforce the tangential boundary conditions strongly, while the normal boundary conditions are enforced weakly through \cref{eq:wrongWeakDiv}. Unfortunately, \cref{eq:wrongDiscreteFormulation} does in general not have a unique solution, as the following counterexample shows. 

Assume $\Omega$ is the unit-square $[0,1]^2$ and let $\mathbf{x}_1, \dots, \mathbf{x}_4$ be its four vertices:
\begin{align*}
    \mathbf{x}_1 &= (0,0), & \mathbf{x}_2 &= (1,0), & \mathbf{x}_3 &= (1,1), & 
    \mathbf{x}_4 &= (0,1).
\end{align*} 
We consider the simplicial mesh made by two triangles $T_1$ and $T_2$.
The triangle $T_1$ has vertices $\mathbf{x}_1, \mathbf{x}_2$ and $\mathbf{x}_4$ while $T_2$ has vertices $\mathbf{x}_2, \mathbf{x}_3$ and $\mathbf{x}_4$. Let $\lambda_1, \dots, \lambda_4$ be the \lq\lq hat functions \rq\rq\, relative to $\mathbf{x}_1, \dots, \mathbf{x}_4$. Let $Q_h$ and $\mathring{V}_h$ be the spaces of Lagrange polynomials of degree one with zero average and the space of Whitney $1$-forms with vanishing tangential components on $\partial \Omega$ respectively. These spaces can be alternatively characterized as 
\begin{align*}
    Q_h &= \operatorname{span}\{ \lambda_1, \dots, \lambda_4\} \cap \HoneAvgZero, \\
    \mathring{\mathbf{V}}_h &= \operatorname{span}\{ \lambda_2 \Grad \lambda_4 - \lambda_4\Grad \lambda_2 \}.
\end{align*}
Set $ \widetilde{q}_h \coloneqq \lambda_1 - \frac16$. Note that $\widetilde{q}_h$ belongs to $Q_h$ since $\int_\Omega (\lambda_1 - \frac16) \,\mathrm{dx} = \frac16 - \frac16 = 0$.
Moreover we have that 
\begin{align*}
    (\Grad \widetilde{q}_h, \lambda_2\Grad \lambda_4 - \lambda_4\Grad \lambda_2) &= \int_{T_1}\Grad \lambda_1 \cdot ( \lambda_2 \Grad \lambda_4 - \lambda_4\Grad \lambda_2)\, \mathrm{dx} 
    \\
    &= \int_{T_1}(- \Grad \lambda_2 - \Grad \lambda_4) \cdot ( \lambda_2 \Grad \lambda_4 - \lambda_4\Grad \lambda_2)\, \mathrm{dx}  \\
    &= \int_{T_1}(\lambda_4 - \lambda_2)\,\mathrm{dx} 
    = 0.
\end{align*}
Thus we have shown that $(\Grad \widetilde{q}_h, \mathbf{v}_h) = 0$ for each $\mathbf{v}_h \in \mathring{\mathbf{V}}_h$. It follows that the non-trivial pair $(0, \widetilde{q}_h)\in\mathring{V}_h \times {Q}_h$ lies in the kernel of the bilinear form in \cref{eq:wrongDiscreteFormulation}. This example can be generalized to arbitrarily fine meshes, and a similar argument holds in three dimensions.
\begin{remark}
    The preceding example necessitates the presence of triangles with two edges on the boundary. This geometric configuration is known to compromise also the well-posedness of the Scott--Vogelius element pair \cite{ScottVogelius}, whose stability relies on delicate mesh conditions. In particular, the inf-sup condition may fail in the presence of such boundary-adjacent triangles.  
\end{remark}

\subsection{Weak imposition of no-slip boundary conditions}
To avoid the problems noted in the previous subsection, we now propose a Nitsche-type method. This leads to the following discrete problem: find $\bfu_h \in \Vh$ and $p_h \in \Qh$ satisfying
\begin{subequations} \label{eq:discreteSystem}
\begin{align*} 
    a_h(\bfu_h, \bfv_h) + b(\bfv_h, p_h) &= l_h(\bfv_h),   & \forall \bfv_h &\in \Vh,\\
    b(\bfu_h, q_h) &= 0, & \forall q_h&\in \Qh.
\end{align*}
\end{subequations}
Here, the bilinear forms $a_h:\Vh\times\Vh\mapsto\R$, $b:\Hcurl\times \HoneAvgZero\mapsto \R$, and functional $l_h:\Vh\mapsto\R$ are defined as
\begin{align}
    a_h(\bfu_h,\bfv_h) &\coloneqq ( \nabla\times \bfu_h, \nabla \times \bfv_h ) + \langle\ntrace(\bfu_h), \ttrace(\bfv_h) \rangle \nonumber\\
    &\quad + \langle\ttrace(\bfu_h), \ntrace(\bfv_h) \rangle+\frac{C_w}{h} \langle\ttrace(\bfu_h), \ttrace(\bfv_h) \rangle, \\
    b(\bfv_h, q_h) &\coloneqq (\bfv_h, \nabla q_h),\\
    l_h(\bfv_h) &\coloneqq (\bff,\bfv_h)+\frac{C_w}{h}\langle  \bfg,\ttrace(\bfv_h)\rangle + \langle  \bfg,\ntrace(\bfv_h)\rangle,
\end{align}
with $C_w>0$ a user-defined constant that needs to be large enough to obtain stability as we will see in the following section. Note that $a_h(\cdot, \cdot)$ is a natural generalization of the bilinear form proposed by Nitsche for the Laplace operator \cite{NitscheBP}. For an extension of this approach to Maxwell's equations, see also \cite[Chapter 45.2]{EG2}.

\subsection{Stability}

To facilitate the stability analysis, we define the following mesh-dependent norm for the velocity:
\begin{equation} \label{eq: hashtagnorm}
    \lVert \bfv_h\rVert^2_{\#} \coloneqq \lVert \bfv_h\rVert^2_{\curl, \Omega} + \frac{1}{h}\lVert \ttrace(\bfv_h)\rVert_\Gamma^2+ h\lVert\ntrace(\bfv_h)\rVert_\Gamma^2,
\end{equation}

Note that $\lVert \cdot \rVert_{\#}$ is well-defined on $\VectorHtwo$ due to \cref{thm:TraceBounded}. In preparation of the error analysis of \Cref{sec: error analysis}, we define the following function space:
\begin{align*}
    \mathbf{V}_{\#} \coloneqq \VectorHtwo + \Vh
\end{align*}
We may extend the domains of the bilinear forms $a_h$ and $b$ to $\mathbf{V}_{\#}$, as shown in the following lemma.

\begin{lemma}[Continuity] \label{lem: continuity}
    The bilinear forms satisfy the following upper bounds:
    \begin{align*}
        a_h(\bfu,\bfv) &\lesssim \| \bfu \|_{\#} \| \bfv \|_{\#}, &
        \forall \bfu,\bfv &\in \mathbf{V}_{\#}, \\
        b( \bfv, q) &\leq\lvert q\rvert_{1,\Omega} \| \bfv \|_\Omega
        \leq \lvert q\rvert_{1,\Omega} \| \bfv \|_{\#}, &
        \forall (q, \bfv) &\in \HoneAvgZero \times \mathbf{V}_{\#}.
    \end{align*}
\end{lemma}
\begin{proof}
    For the continuity of $a_h$, we apply Cauchy-Schwarz to each term and introduce a scaling with $h$ where necessary. 
    If we subsequently apply Cauchy-Schwarz to the sum of products, we obtain
    \begin{align*}
        a_h(\bfu,\bfv)
        &\leq \| \nabla \times \bfu \|_{\Omega} \| \nabla \times \bfv \|_{\Omega} + 
        h^{-\frac{1}{2}} \|\ntrace( \bfu )\|_\Gamma h^{\frac{1}{2}} \|\ttrace( \bfv )\|_\Gamma \nonumber\\
        & \quad
        + h^{-\frac{1}{2}} \|\ttrace( \bfu )\|_\Gamma h^{\frac{1}{2}} \|\ntrace( \bfv )\|_\Gamma 
        + C_w h^{-\frac{1}{2}} \|\ttrace( \bfu )\|_\Gamma h^{-\frac{1}{2}} \|\ttrace( \bfv )\|_\Gamma
        \nonumber\\
        &\lesssim 
        \| \bfu \|_{\#} \| \bfv \|_{\#}
    \end{align*}
    
    The continuity of $b$ follows by a similar application of Cauchy-Schwarz.
\end{proof}

\begin{theorem}[Coercivity] \label{thm:coercivity}
For $C_w>0$ big enough and $h>0$ small enough, there exists a $C>0$ such that for all $\bfu_h \in \X_h$
\begin{equation} \label{eq:biggerzero}
    a_h(\bfu_h,\bfu_h) \geq C \lVert \bfu_h\rVert_{\#}^2.
\end{equation}

\end{theorem}
\begin{proof}
    Assuming $\bfu_h \in \X_h$ given, we proceed in three steps.
\begin{itemize}[leftmargin=*]
    \item Step 1: A lower bound on $a_h(\bfu_h,\bfu_h)$. We use Young's inequality and the discrete trace inequality from \Cref{lem: discrete trace ineq} to obtain
    \begin{align*}
        \langle \ntrace(\bfu_h), \ttrace(\bfu_h) \rangle 
        &\geq - \|\ntrace(\bfu_h)\|_{\Gamma}\|\ttrace(\bfu_h)\|_{\Gamma} \\
        &\geq - \frac12 \left(\frac{1}{\alpha}\|\ntrace(\bfu_h)\|_\Gamma^2 + \alpha \|\ttrace(\bfu_h)\|_\Gamma^2 \right) \\
        &\geq - \frac12 \left(\frac{C_n}{\alpha h}\|\nabla\times\bfu_h\|^2_{\Omega} + \alpha \|\ttrace(\bfu_h)\|_\Gamma^2 \right),
    \end{align*}
    where $\alpha >0$ is a constant to be chosen later. We then derive
    \begin{align*} 
        a_h(\bfu_h,\bfu_h) &= \| \nabla\times\bfu_h \|_{\Omega}^2 + 2\langle \ntrace(\bfu_h), \ttrace(\bfu_h) \rangle + \frac{C_w}{h} \| \ttrace(\bfu_h) \|_\Gamma^2 \\
        &\geq \| \nabla\times\bfu_h \|_{\Omega}^2 - \left(\frac{C_n}{\alpha h}\|\nabla\times\bfu_h\|^2_{\Omega} + \alpha\|\ttrace(\bfu_h)\|_\Gamma^2 \right) + \frac{C_w}{h} \| \ttrace(\bfu_h) \|_\Gamma^2 \\
        &=  \left( 1-\frac{C_n}{\alpha h}\right) \|\nabla\times\bfu_h\|^2_{\Omega} + \left(\frac{C_w}{h}-\alpha\right)\|\ttrace(\bfu_h)\|_\Gamma^2.
    \end{align*}
    We choose $\alpha = 2\frac{C_n}{h}$ and $C_w \ge 2\alpha h = 4C_n$. Then, 
    \begin{align} \label{eq:lowerbounda}
        a_h(\bfu_h,\bfu_h) &\geq  \frac{1}{2} \|\nabla\times\bfu_h\|^2_{\Omega} + 2\frac{C_n}{h}\|\ttrace(\bfu_h)\|_\Gamma^2.
    \end{align}
    
    \item Step 2: An upper bound on $\|\bfu_h\|^2_{\Omega}$. Since $\mathbf{u}_h\in \X_h$, using \eqref{eq:HH_decomp} there exist $\bfu_h^{\perp}\in \mathbf{Z}_h$ and a discrete harmonic 1-form $\mathfrak{h}_h \in \mathfrak{H}^1_h$ such that
    \begin{equation*}
        \bfu_h = \bfu_h^\perp + \mathfrak{h}_h.
    \end{equation*}
    By the discrete Poincaré inequality\cite[Thm.~4.7]{Hip02} there exists a constant $C_P>0$ such that
    \begin{equation}\label{eq:bounduperp}
        \|\bfu_h^\perp\|_{\Omega} \leq C_P \|\nabla\times\bfu_h^\perp\|_{\Omega}.
    \end{equation}
       Combining \eqref{eq:bounduperp}, \cref{lemma:boundHmh}, and \cref{lem:L2BoundaryBoundDiscrete}, we obtain
    \begin{align*}
        \|\bfu_h\|_{\Omega} &\leq \|\bfu_h^\perp\|_{\Omega} 
        + \|\mathfrak{h}_h\|_{\Omega}\\
        &\lesssim  \|\bfu_h^\perp\|_{\Omega} 
        + \|\ttrace(\mathfrak{h}_h)\|_{-\frac{1}{2},t, \Gamma}\\
        &\leq  \| \bfu_h^\perp\|_{\Omega} 
        + \|\ttrace(\bfu_h)\|_{-\frac{1}{2},t,\Gamma}+\|\ttrace(\bfu_h^\perp)\|_{-\frac{1}{2}, \Gamma}\\
        &\lesssim \|\bfu_h^\perp\|_{\curl, \Omega}
        + \|\ttrace(\bfu_h)\|_{\Gamma} \\
        &\lesssim  \|\nabla\times\bfu_h^\perp\|_{\Omega} 
        + \|\ttrace(\bfu_h)\|_{\Gamma}\\
        &=  \|\nabla\times\bfu_h\|_{\Omega} 
        + \|\ttrace(\bfu_h)\|_{\Gamma},
    \end{align*}

    \item Step 3: By combining step 1 and 2 with the discrete trace inequality \eqref{eq:trace_ineq_gamman}, we obtain 
    \begin{align*}
        \lVert \bfu_h\rVert_{\#}^2
        &= \|\bfu_h\|_{\curl, \Omega}^2 + \frac{1}{h}\|\ttrace(\bfu_h)\|_\Gamma^2 + h\|\ntrace(\bfu_h)\|_\Gamma^2\\
        &\lesssim \|\nabla \times \bfu_h\|_\Omega^2 + \frac{1}{h}\|\ttrace(\bfu_h)\|_\Gamma^2 \\
        &\lesssim a_h(\bfu_h,\bfu_h).
    \end{align*}
\end{itemize}
    This completes the proof.
\end{proof}
\begin{lemma}[Inf-sup] \label{lem: infsup b}
    There exists $\beta>0$ independent of $h$ such that 
\begin{equation*}
    \inf_{q_h \in \Qh} \sup_{\bfv_h \in \Vh} 
    \frac{ b(\bfv_h, q_h)}{\lvert q_h \rvert_{1,\Omega} \lVert \bfv_h \rVert_{\#} }
    \geq \beta h > 0.
\end{equation*}
\end{lemma}
\begin{proof}
Let $q_h \in \Qh$. Since $\nabla q_h \in \Vh$, we may apply the discrete trace inequality\eqref{eq:trace_ineq_gammapar} to bound
\begin{equation*}
\begin{split}
    \lVert \Grad q_h \rVert_{\#}^2 & = \lVert \Grad q_h \rVert_{\Omega}^2 + \frac{1}{h}\lVert \ttrace(\Grad q_h) \rVert_\Gamma^2 \\ 
    & \leq \lVert \Grad q_h \rVert_{\Omega}^2 + C_\parallel\frac1{h^2} \lVert \Grad q_h \rVert_\Omega^2
    \lesssim h^{-2} \lVert \Grad q_h \rVert_\Omega^2
    = h^{-2} \lvert q_h \rvert_{1,\Omega}^2.
    \end{split}
\end{equation*}
Taking $\bfv_h \coloneqq\Grad q_h$ we obtain
\begin{align*}
    \sup_{\bfv_h \in \Vh} \frac{ b( \bfv_h, q_h)}{\lVert \bfv_h \rVert_{\#}} &\geq \frac{ b( \Grad q_h,q_h)}{\lVert \Grad q_h \rVert_{\#}}
    = \frac{ \lVert \Grad q_h \rVert_{\Omega}^2}{\lVert \Grad q_h \rVert_{\#}} 
    = \frac{ \lvert q_h\rvert_{1,\Omega}^2}{\lVert \Grad q_h \rVert_{\#}} 
    \gtrsim h \lvert q_h\rvert_{1,\Omega},
\end{align*}
concluding the proof.
\end{proof}

\begin{theorem}[Stability]
    The discrete problem \eqref{eq:discreteSystem} admits a unique solution $(\bfu_h, p_h) \in \Vh \times \Qh$ that satisfies
    \begin{align*}
        \| \bfu_h \|_{\#} + h \| \nabla p_h \|_\Omega 
        &\lesssim \| l_h \|_{\Vh'} \coloneqq \sup_{\bfv \in \Vh \setminus \{\boldsymbol{0}\}}\frac{ \lvert \ell_h(\bfv_h) \rvert}{\lVert \bfv_h \rVert_{\#}}  
    \end{align*}
\label{thm: stability}
\end{theorem}
\begin{proof}
    \Cref{lem: continuity,thm:coercivity,lem: infsup b} suffice to invoke saddle point theory. The $h$-dependency of the inf-sup constant in \cref{lem: infsup b} is reflected by the scaling on the pressure term \cite[Thm.~4.2.3]{boffi2013mixed}.
\end{proof}

\begin{remark}
As we will see in \cref{lem: suboptimal p conv}, the presence of the $h$ factor in \cref{lem: infsup b} and \cref{thm: stability} implies that the error of the pressure converges slower than the best approximation error in $H^1$. 
\end{remark} 

\subsection{A priori error analysis}
\label{sec: error analysis}

We proceed by deriving consistency and error estimates. Recall that \Cref{lem: continuity} showed the continuity of $a_h$ with respect to the $\lVert \cdot \rVert_{\#}$ norm from \eqref{eq: hashtagnorm}. We are now ready to state the following consistency estimate.
\begin{lemma}[Consistency] \label{thm:Consistency}
    Let $(\bfu,p)\in \VectorHtwo\times H^1(\Omega)$ be the solution of \eqref{eq:continuousSymGradStokes}. Then, we have consistency in the following sense:
    \begin{align*}
        a_h(\bfu,\bfv_h)+b(\bfv_h,p) &= l_h(\bfv_h),&\forall \bfv_h \in \Vh,\\
        b( \bfu, q_h) &=0,&\forall q_h \in \Qh.
    \end{align*}
\end{lemma}
\begin{proof}
    Since $\bfu\in \VectorHtwo$ solves \eqref{eq:continuousSymGradStokes}, we have for any $\bfv_h \in \Vh$:
    \begin{align*} 
        l_h(\bfv_h) &= (\bff,\bfv_h) + \langle  \bfg,\ntrace(\bfv_h)) + \frac{C_w}{h}\langle \bfg,\ttrace(\bfv_h)\rangle \\
        &=  (\nabla\times\nabla\times\bfu+\nabla p,\bfv_h) + \langle \bfg,\ntrace(\bfv_h)\rangle + \frac{C_w}{h}\langle \bfg,\ttrace(\bfv_h)\rangle\\
        &=  ( \nabla\times\bfu, \nabla\times\bfv_h ) +( \nabla p,\bfv_h) +  \langle\ntrace(\bfu),\bfv_h\rangle 
        + \langle\bfg,\ntrace(\bfv_h)\rangle + \frac{C_w}{h}\langle \bfg,\ttrace(\bfv_h)\rangle\\
        &= a_h(\bfu, \bfv_h) + b(\bfv_h, p),
    \end{align*}
    where we used integration by parts and the fact that $\bfu=\bfg$ on $\Gamma$.
\end{proof}

For ease of reference, we state the following identity.
\begin{corollary}[Orthogonality] \label{cor: orthogonality}
    The continuous and discrete solutions satisfy
    \begin{align*}
        a_h(\bfu - \bfu_h,\bfv_h)+b(\bfv_h, p - p_h) &= 0,
        &\forall \bfv_h \in \Vh.
    \end{align*}
\end{corollary}

\subsubsection{A priori error estimates in the natural norm}

\begin{theorem} \label{thm:convhashtagnorm}
    Let $(\bfu,p)\in \VectorHtwo \times \HoneAvgZero$ solve \eqref{eq:continuousSymGradStokes}. Then, if $(\bfu_h,p_h)\in \Vh\times \Qh$ are Galerkin solutions of \eqref{eq:discreteSystem}, we have 
    \begin{equation} \label{eq:optimalityInfimum}
        \|\bfu-\bfu_h\|_{\#} \lesssim \inf_{\bfv_h \in \X_h}\|\bfu-\bfv_h\|_{\#} + \inf_{q_h \in \Qh}\lvert p-q_h\rvert_{1,\Omega}.
    \end{equation}
\end{theorem}

\begin{proof}
    Let $\bfv_h \in \X_h$, then, using \cref{thm:coercivity} and \cref{cor: orthogonality}, we can estimate
    \begin{align*}
        \|\bfu_h-\bfv_h\|_{\#} &\lesssim \sup_{\bfw_h \in \X_h} \frac{a_h(\bfu_h-\bfv_h,\bfw_h)}{\|\bfw_h\|_{\#}}
        = \sup_{\bfw_h \in \X_h} \frac{b(\bfw_h, p - p_h)+a_h(\bfu-\bfv_h,\bfw_h)}{\|\bfw_h\|_{\#}},
    \end{align*}
    Since $\bfw_h \in \X_h$, we have for all $q_h \in \Qh$
    \begin{align*}
        b(\bfw_h,p - p_h,) = b(\bfw_h, p-q_h) \leq \lvert p-q_h\rvert_{1,\Omega}\|\bfw_h\|_{\#},
    \end{align*}
    where the inequality is the continuity of $b$ from \Cref{lem: continuity}.    
    Similarly, the continuity of $a_h$ gives us
    \begin{align*}
        a_h(\bfu-\bfv_h,\bfw_h) 
        &\lesssim \|\bfu-\bfv_h\|_{\#}\|\bfw_h\|_{\#}.
    \end{align*}
    We thus find
    \begin{align*}
        \|\bfu_h-\bfv_h\|_{\#} \lesssim\|\bfu-\bfv_h\|_{\#} + \lvert p-q_h\rvert_{1,\Omega}
    \end{align*}
    and we conclude that for all $\bfv_h \in \X_h$ and $q_h \in \widetilde{H}^1(\Omega)$
    \begin{align*}
        \|\bfu-\bfu_h\|_{\#} \leq \|\bfu-\bfv_h\|_{\#}+\|\bfv_h-\bfu_h\|_{\#}\lesssim(\|\bfu-\bfv_h\|_{\#} + \lvert p-q_h\rvert_{1,\Omega}).
    \end{align*}
    This completes the proof.
\end{proof}
\begin{lemma} \label{lem: suboptimal p conv}
    The following error estimate holds:
    \begin{equation*}
        \lVert \Grad (p - p_h) \rVert_{\Omega} \lesssim 
        h^{-1}\left(\inf_{\bfv_h \in \X_h } \lVert \bfu-  \bfv_h \rVert_{\#} + \inf_{q_h \in \Qh}\lvert p - q_h \rvert_{1,\Omega}\right).
    \end{equation*}
\end{lemma}
\begin{proof}
Using the inf-sup condition from \cref{lem: infsup b}, the consistency from \cref{cor: orthogonality}, and the continuity from \cref{lem: continuity}, we get
\begin{align*}
\lVert \Grad ( p_h - q_h )\rVert_{\Omega} & \leq \frac{1}{\beta h}\sup_{\bfv_h \in \Vh} \frac{b( \bfv_h, p_h - q_h)}{\lVert \bfv_h \rVert_{\#}} \\
&= \frac{1}{\beta h}\sup_{\bfv_h \in \Vh}  \frac{ a_h(\bfu - \bfu_h, \bfv_h) + b(\bfv_h, p-q_h)}{\lVert \bfv_h \rVert_{\#} } \\ 
		& \lesssim h^{-1} \left( \| \bfu - \bfu_h \|_{\#} +\lvert p - q_h \rvert_{1,\Omega} \right) 
\end{align*}
Combining this bound with \cref{thm:convhashtagnorm} and a triangle inequality, we obtain the result. 
\end{proof}

Thus, we get an a priori estimate in the pressure that is one entire order less with respect to the velocity, due to the dependency of the inf-sup constant on $h$ in \Cref{lem: infsup b}. 

The above results involve the best approximation error on $\X_h$, but in practice it is more useful to have an estimate which involves the best approximation error on $\Vh$.  To fix the ideas, we assume now that $\Omega_h$ is simplicial,  $\Qh$ is the space of Lagrange polynomials of degree $r$ and $\Vh$ is the space of N\'{e}d\'{e}lec of the first kind \cite{Ned80} of degree $r$ (i.e. Whitney elements correspond to $r = 1$). Let $\mathcal{I}^c_h$ be the standard commuting interpolation operator onto $\Vh$. We summarize the local approximation properties of $\mathcal{I}^c_h$ in the following Theorem. 
\begin{theorem}
Let $q\in[2, \infty].$ Set $Z(T) \coloneqq  \mathbf{W}^{q,s}(T)$ with $qs>1$ if $d = 2$ and $qs>2$ if $d = 3$. Then $\mathcal{I}^c_h:Z(T) \to \Vh(T)$ satisfies the local approximation estimates 
\begin{subequations}
    \begin{align}
        \lvert \bfv - \mathcal{I}^c_h\bfv\rvert_{m, T, q} &\lesssim \begin{cases} h^{r-m}_T \lvert \bfv \rvert_{r,T,q} \text{ if $r > 1$ or $q>2$ or $d = 2$,}\\ 
        h^{1-m}_T\lvert \bfv \rvert_{1,T,q} + h^{2-m}_T\lvert \bfv \rvert_{2,T,q} \text{ if $r = 1$ and $q = 2$ and $d = 3$,}\end{cases} \label{eq:local_approx_1}\\
        \lvert \Curl (\bfv - \mathcal{I}^c_h\bfv)\rvert_{m,T, q} &\lesssim h^{r-m}_T\lvert \Curl \bfv \rvert_{r,T,q}, \label{eq:local_approx_2}
    \end{align}
\end{subequations}
for each $T\in\Omega_h$. Here $m\in \{0,\dots, r\}$ is an integer.
\end{theorem} 
\begin{proof} 
See \cite[Theorem 16.10]{EG1}.
\end{proof}
We can use property \eqref{eq:local_approx_1} in combination with the multiplicative trace inequality \cite[Lemma 12.15]{EG1} to obtain the following approximation property on faces. 
\begin{lemma} Let $F$ be a face of the element $T$, and $\bfv \in \mathbf{W}^{r, q}(T)$. Assume $q>2$ or $r>1$ or $d = 2$.  Then 
\begin{equation}
    \lVert \ttrace( \bfv - \mathcal{I}^c \bfv ) \rVert_{F,q} \lesssim h^{r- \frac{1}{q}}_T \lvert \bfv \rvert_{r,T,q}.
    \label{eq:local_trace_approx}
\end{equation}
\end{lemma}
\begin{proof}We have
\begin{align*}
    \lVert \ttrace(\bfv - \mathcal{I}^c_h\bfv ) \rVert_{F,q} & \lesssim h^{-\frac{1}{q}}_T\lVert \bfv - \mathcal{I}^c_h\bfv \rVert_{\mathbf{L}^q(T)} + \lVert \bfv - \mathcal{I}^c_h\bfv  \rVert_{\mathbf{L}^q(T)}^{1-\frac{1}{q}}\lvert \bfv - \mathcal{I}^c_h\bfv  \rvert_{1,T,q}^{\frac{1}{q}}\\
    & \lesssim h^{r-\frac{1}{q}}_T \lvert \bfv \rvert_{r,T,q} + h_T^{r\left(1 - \frac{1}{q}\right)} h_T^{\frac{r-1}{q}}\lvert \bfv \rvert_{r,T,q} \\
    &=  h^{r-\frac{1}{q}}_T \lvert \bfv \rvert_{r,T,q} + h_T^{r\left(1 - \frac{1}{q}\right)} h_T^{\frac{r-1}{q}}\lvert \bfv \rvert_{r,T,q} \\
    &= 2 h^{r-\frac{1}{q}}_T \lvert \bfv \rvert_{r,T,q}.
\end{align*}
\end{proof}
The case $q = 2$, $r = 1$ and $d = 3$ is similar. We omit it for brevity. 
Since the mesh is assumed to be quasi-uniform, we have also the global approximation estimates
\begin{subequations}
\begin{align}
    \lVert \bfv - \mathcal{I}^c_h \bfv\rVert_{\Omega,q} &\lesssim h^r \lvert \bfv \rvert_{r,\Omega,q}, \\ 
    \lVert \Curl\bfv -  \Curl \mathcal{I}^c_h\bfv\rVert_{\Omega,q} &\lesssim h^r \lvert \Curl \bfv\rvert_{r,\Omega,q},
\end{align}
\label{eq:global_approx}
\end{subequations}
Moreover, \eqref{eq:local_trace_approx} implies the global bound
\begin{equation}
    \lVert \ttrace(\bfv - \mathcal{I}^c_h\bfv )\rVert_{\Gamma,q} \lesssim h^{r-\frac{1}{q}} \lvert \bfv \rvert_{r,\Omega,q}.
    \label{eq:global_trace_approx}
\end{equation}
Similarly, using \eqref{eq:local_approx_2}, it is possible to show 
\begin{equation}
    \lVert \ntrace(\bfv -  \mathcal{I}^c_h \bfv )\rVert_{\Gamma,q} \lesssim h^{r- \frac{1}{q}}\lvert \Curl \bfv \rvert_{r,\Omega,q}.
    \label{eq:global_trace_approx_gamman}
\end{equation} 
With these technical results, we are ready to state and prove the a priori error estimates of our method.

 \begin{theorem} \label{thm: conv hashtag norm}
     Let $(\mathbf{u}, p)$ be the exact solution of \eqref{eq:continuousSymGradStokes}, and assume $\mathbf{u}\in \mathbf{H}^{r}(\Omega)$, $\Curl \bfu \in \mathbf{H}^r(\Omega)$ and $p\in H^{r+1}(\Omega)$. Then the following error estimate holds:
     \begin{equation} \label{eq: conv hashtag nd}
         h \lVert \Grad (p - p_h) \Vert_{\Omega} + \lVert \bfu - \bfu_h \rVert_{\#} \lesssim h^{r-1}\lvert \bfu \rvert_{r, \Omega} + h^{r}\lvert \Curl \bfu \rvert_{r, \Omega} + h^r\lvert p \rvert_{r+1, \Omega}. 
     \end{equation}
     For two-dimensional, convex domains, and assuming $\bfu \in \mathbf{W}^{r, \infty}(\Omega)$ this estimate can be improved to 
     \begin{equation} \label{eq: conv hashtag 2d}
         h \lVert \Grad (p - p_h) \Vert_{\Omega} + \lVert \bfu - \bfu_h \rVert_{\#} \lesssim h^{r-\frac{1}{2}} \lvert \bfu \rvert_{r,\Omega,\infty} +h^r\lvert \Curl \bfu \rvert_{r, \Omega} + h^r \lvert p \rvert_{r+1, \Omega}. 
     \end{equation}
\end{theorem}
\begin{proof}
    Let $P_{\X_h}$ and $P_{\Grad \Qh}$ be the orthogonal projectors onto $\X_h$ and $\Grad \Qh$ respectively. In particular, note that $\mathbf{v}_h = P_{\X_h}\mathbf{v}_h + P_{\Grad \Qh}\mathbf{v}_h$ for each $\mathbf{v}_h \in \Vh$. Then we have
\begin{equation*}
\begin{split}
    \inf_{\bfz_h \in \X_h}\lVert \bfu - \bfz_h\rVert_{\#} & = \inf_{\bfv\in \Vh}\lVert \bfu - P_{\X_h}\bfv_h\rVert_{\#}\\
    &= \inf_{\bfv_h \in \Vh}\lVert \bfu - \bfv_h\rVert_{\#} + \lVert P_{\Grad \Qh}\bfv_h\rVert_{\#}\\
    &= \inf_{\bfv_h \in \Vh}\lVert \bfu - \bfv_h\rVert_{\#} + \lVert P_{\Grad \Qh}(\bfu -\bfv_h)\rVert_{\#} \\
    &\eqqcolon\inf_{\bfv_h \in \Vh} \mathrm{I} + \mathrm{II}.
    \end{split}
\end{equation*}

We take $\bfv_h = \mathcal{I}^c_h\bfu$. Then, we can bound the boundary terms in $\mathrm{I}$ using the embedding of $L^q(\Gamma)$ in $L^2(\Gamma)$ for each $q\in [2, \infty]$, and the bounds \eqref{eq:global_trace_approx} and \eqref{eq:global_trace_approx_gamman}:
\begin{align*}
\lVert \ttrace(\bfu - \mathcal{I}^c_h\bfu) \rVert_{\Gamma}&\leq \lvert \Gamma \rvert^{\frac{1}{2} - \frac{1}{q}} \lVert \ttrace(\bfu - \mathcal{I}^c_h\bfu) \rVert_{\Gamma,q} \\ 
& \lesssim \lvert \Gamma \rvert^{\frac{1}{2} - \frac{1}{q}} h^{r- \frac{1}{q}} \lvert \bfu \rvert_{r,\Omega,q};  \\
\lVert \ntrace(\bfu - \mathcal{I}^c_h\bfu)\rVert_{\Gamma}& \leq \lvert \Gamma \rvert^{\frac{1}{2} - \frac{1}{q}} \lVert \ntrace(\bfu - \mathcal{I}^c_h\bfu) \rVert_{\Gamma,q} \\ 
& \lesssim  \lvert \Gamma \rvert^{\frac{1}{2} - \frac{1}{q}} h^{r- \frac{1}{q}} \lvert \Curl \bfu \rvert_{r,\Omega,q},
\end{align*}
provided that the norms on the right-hand side are finite. Then, the interior terms in $\mathrm{I}$ can be bounded with \eqref{eq:global_approx}. It follows that
\begin{equation*}
    \mathrm{I} \lesssim h^r\lvert \bfu \rvert_{r, \Omega} + h^r\lvert \Curl \bfu\rvert_{r, \Omega} + h^{r + \frac{1}{2} - \frac{1}{q}}\lvert \Curl \bfu\rvert_{r,\Omega,q} + h^{r- \frac{1}{2}- \frac{1}{q}}\lvert \bfu \rvert_{r,\Omega,q}.
\end{equation*}
We now estimate $\mathrm{II}$. The only nonzero terms are 
\begin{align*}
    \lVert P_{\Grad \Qh}(\mathbf{u} - \mathcal{I}^c_h \mathbf{u}) \rVert_{\Omega} &\leq \lVert \mathbf{u} - \mathcal{I}^c_h \mathbf{u} \rVert_{\Omega}
    \lesssim h^r \lvert \bfu \rvert_{r, \Omega}, \\ 
    \lVert\ttrace P_{\Grad \Qh}(\mathbf{u} - \mathcal{I}^c_h \mathbf{u}) \rVert_{\Gamma} &\lesssim h^{-\frac{1}{2}} \lVert P_{\Grad \Qh}(\mathbf{u} - \mathcal{I}^c_h \mathbf{u}) \rVert_{\Omega}
    \lesssim h^{r-\frac{1}{2}} \lvert \bfu \rvert_{r, \Omega}.
\end{align*}
It follows that 
\begin{equation*}
    \mathrm{II} \lesssim h^r \lvert \bfu \rvert_{r, \Omega} + h^{r-1}\lvert \bfu \rvert_{r, \Omega}.
\end{equation*}
 Taking $q = 2$ and $q_h$ in \cref{lem: suboptimal p conv} as the Lagrange interpolant of $p$, we get the estimate \eqref{eq: conv hashtag nd}.
 In dimension two, when $\Omega$ is convex and $\bfu \in\mathbf{W}^{r, \infty}(\Omega)$, the last bound can be improved thanks to \cite[Corollary 3.1]{Duran}:
 \begin{align*}
     \lVert \ttrace P_{\Grad \Qh}(\mathbf{u} - \mathcal{I}^c_h \mathbf{u}) \rVert_{\Gamma} & \leq \lvert \Gamma \rvert^{\frac{1}{2} - \frac{1}{q}}\lVert P_{\Grad \Qh}(\mathbf{u} - \mathcal{I}^c_h \mathbf{u}) \rVert_{\Gamma,q}\\
     & \lesssim\lvert \Gamma \rvert^{\frac{1}{2} - \frac{1}{q}} h^{-\frac{1}{q}}\lVert P_{\Grad \Qh}(\mathbf{u} - \mathcal{I}^c_h \mathbf{u}) \rVert_{\Omega,q} \\
     &  \lesssim\lvert \Gamma \rvert^{\frac{1}{2} - \frac{1}{q}} h^{-\frac{1}{q}}q \lVert \bfu - \mathcal{I}^c_h \bfu\rVert_{\Omega,q} \\ 
     & \lesssim \lvert \Gamma \rvert^{\frac{1}{2} - \frac{1}{q}} h^{r-\frac{1}{q}}q\lVert \bfu \rVert_{r,\Omega,q}.
 \end{align*}
 Taking $q = \lvert \log h\rvert$, we obtain \eqref{eq: conv hashtag 2d}.
\end{proof}
\begin{remark} \label{rem:improvedconvergence3D}
    The improved error estimates of \Cref{thm: conv hashtag norm} rely on the $L^q$-stability of the projection $P_{\Grad Q_h}$. To the best of the authors' knowledge, this result has been established in the literature only for two dimensional topologically trivial domains that are either convex or have smooth boundary. We conjecture that the stability property extends to three dimensions and to non-topologically trivial domains, as our numerical experiments do not reveal any deterioration in the observed convergence rates; see Section~\ref{sec:numerics}.
\end{remark}

\subsubsection{Improved estimates in \texorpdfstring{$L^2$}{L2} using duality techniques}
In this section we derive improved $L^2$ error estimates for pressure and velocity.  We need the averaging interpolator with boundary prescription $\iEG: \mathbf{L}^1(\Omega) \to \mathring{\mathbf{V}}_h$ introduced and analyzed by Ern and Guermond \cite{EG17}. In particular, we recall that $\iEG$ is stable in $L^2$ and satisfies the approximation estimate
\begin{equation}
\lVert \bfv - \iEG(\bfv)\rVert_{\Omega} \lesssim h\lVert \bfv \rVert_{1, \Omega},
\label{eq:approx_EG}
\end{equation}
for any $\bfv \in \mathbf{H}^1_0(\Omega)$.

\begin{theorem} \label{lem: p conv L2}
    The error in the pressure satisfies the following bound
    \begin{equation*}
        \lVert p - p_h \rVert_{\Omega} \lesssim \lVert \bfu - \bfu_h \rVert_{\#} + \inf_{q_h \in \Qh}\lvert p - q_h \rvert_{1,\Omega}.
    \end{equation*}
\end{theorem}
\begin{proof}
    Using \cite{Bogovskii}, we may construct $\mathbf{w}\in \VectorH^1_0(\Omega)$ such that
    \begin{align} \label{eq: Stokes test func}
        -\Div \mathbf{w} &= p - p_h, &
        \lVert \mathbf{w} \rVert_{1, \Omega} &\lesssim \lVert p - p_h \rVert_{\Omega}.
    \end{align}

    Using the Galerkin orthogonality from \Cref{cor: orthogonality}, we derive
    \begin{align}
        \lVert p - p_h \rVert_{\Omega}^2 &= -(\Div \mathbf{w}, p - p_h)  \nonumber\\
        &= (\mathbf{w}, \Grad (p - p_h)) \nonumber\\ 
        &= (\mathbf{w} - \iEG(\mathbf{w}), \Grad (p - p_h)) + b( \iEG(\mathbf{w}),p - p_h) \nonumber\\
        &= (\mathbf{w} - \iEG(\mathbf{w}), \Grad (p - p_h)) - a_h(\bfu - \bfu_h, \iEG(\mathbf{w})) .\label{eq: intermediate bound}
    \end{align}
    To bound the final term, we refine the continuity result from \Cref{lem: continuity} to the special case of $\bfu \in \mathbf{V}_{\#}$ and $\mathring{\bfv}_h \in \Vhz$. Similar to \Cref{lem: continuity}, we use the Cauchy-Schwarz inequality and the discrete trace inequality \eqref{eq:trace_ineq_gamman} to obtain
    \begin{align*}
        a_h(\bfu, \mathring{\bfv}_h) &= (\Curl \bfu, \Curl \mathring{\bfv}_h) + \langle \ttrace(\bfu), \ntrace(\mathring{\bfv}_h) \rangle \\
        & \leq \lVert \Curl \bfu\rVert_{\Omega}\lVert \Curl \mathring{\bfv}_h \rVert_{\Omega}
        + h^{-\frac{1}{2}}\lVert \ttrace(\bfu) \rVert_\Gamma h^{\frac{1}{2}}\lVert \ntrace(\mathring{\bfv}_h)\rVert_\Gamma \\
        & \lesssim \lVert \bfu \rVert_{\#} \lVert \Curl \mathring{\bfv}_h \rVert_{\Omega}.
    \end{align*}

    We now use this bound and apply \eqref{eq:approx_EG} in order to deduce
    \begin{align*}
        \lVert p - p_h \rVert_{\Omega}^2 
        &\lesssim \| \mathbf{w} - \iEG(\mathbf{w}) \|_\Omega \| \Grad (p - p_h) \|_\Omega + \lVert  \bfu - \bfu_h\rVert_{\#} \lVert \Curl  \iEG(\mathbf{w})\rVert_{\Omega}\\
        &\lesssim h \lVert \mathbf{w} \rVert_{1, \Omega} \lVert \Grad (p - p_h )\rVert_{\Omega} + \lVert \bfu - \bfu_h \rVert_{\#} \lVert \mathbf{w} \rVert_{1, \Omega} \\ 
        & \lesssim \lVert p - p_h \rVert_{\Omega}\left(\lVert \bfu - \bfu_h \rVert_{\#} + \inf_{q_h \in \Qh}\lvert p - q_h \rvert_{1,\Omega} \right).
    \end{align*}
    The final inequality is due to \eqref{eq: Stokes test func} and \Cref{lem: suboptimal p conv}. This implies the result.
\end{proof}

For the velocity, we confine ourselves to the case of a two-dimensional convex domain. Let $\mathbb{P}_r$ be the space of discontinuous, elementwise polynomials of degree $r$.
We follow closely \cite[Thm.~3.7]{ArFaGo2012}, and define a particular projection operator in the following lemma. 

\begin{lemma} \label{lem: projection Vzero}
    Let $P_{\Vhz}: \VectorH_0(\operatorname{curl}, \Omega) \to \Vhz$ be such that, for given $\bfv \in \VectorH_0(\operatorname{curl}, \Omega)$, the projection $P_{\Vhz} \bfv \in \Vhz$ satisfies
    \begin{align*}
        (P_{\Vhz}\bfv, \Grad \mathring{q}_h)  
        &= (\bfv, \Grad \mathring{q}_h), \qquad &\forall \mathring{q}_h \in \mathring{Q}_h,\\
        (\Curl P_{\Vhz} \bfv,  s_h)& = (\Curl \bfv, s_h), \qquad &\forall s_h \in \mathbb{P}_{r-1}. 
    \end{align*}
    Then the following error estimates hold
    \begin{align}
        \lVert \bfv - P_{\Vhz}\bfv \rVert_{\Omega,q} 
        &\lesssim q h^l \lVert \bfv \rVert_{l, \Omega, q}, &
        1\leq l\leq r,\, 2\leq q<\infty, \label{eq:thm3.4a} \\
        \lVert \Curl (\bfv - P_{\Vhz} \bfv )\rVert_{\Omega} 
        &\lesssim h^l \lVert \Curl \bfv \rVert_{l, \Omega},  &
        0\leq l \leq r.\label{eq:thm3.4b}
    \end{align}

    Moreover, for $2\leq q\leq \infty$ and each $q_h \in \Qh$ and $\bfv\in \VectorH_0(\operatorname{curl}, \Omega)\cap \mathbf{L}^q(\Omega)$ it holds 
    \begin{equation}
        (\Grad q_h, \bfv - P_{\Vhz}\bfv) \leq Ch^{-\frac{1}{2} - \frac{1}{q}}\lVert q_h \rVert_{\Omega} \lVert \bfv - P_{\Vhz}\bfv\rVert_{\Omega,q}. \label{eq:thm3.5}
    \end{equation}
\end{lemma}
\begin{proof}
See \cite[Thm.~3.4 and Thm.~3.5]{ArFaGo2012}.
\end{proof}

We continue by introducing a projection operator for the pressure variable.
\begin{lemma} \label{lem: projection Qh}
    Let $P_{\Qh}: H^1(\Omega) \to \Qh$ be the elliptic projection that, for given $q$, solves
    \begin{align*}
        (\Grad P_{\Qh}q, \Grad \varphi_h) &= (\Grad q, \Grad \varphi_h), 
        &\forall \varphi_h &\in \Qh.
    \end{align*}
    Then the following error estimate holds
    \begin{equation*}
        \lVert q - P_{\Qh}q\rVert_{\Omega} + h \lVert \Grad( p - P_{\Qh}p) \rVert_{\Omega} \leq Ch^l \lVert p \rVert_{l, \Omega}, \quad 1 \leq l \leq r.
    \end{equation*}
    Moreover, the following additional property holds for $q\in H^1(\Omega)$ and $\bfv_h \in \Vh$:
    \begin{equation}
        (\Grad(q - P_{\Qh}q), \bfv_h)\leq Ch \lVert \Grad (q - P_{\Qh}q)\rVert_{\Omega}\lVert \Curl \bfv_h \rVert_{\Omega}.
        \label{eq:3.15}
    \end{equation}
\end{lemma}
\begin{proof}
    The error bound is a classical finite element result for the Laplace problem. For a proof of \eqref{eq:3.15}, see \cite[Equation (3.15)]{ArFaGo2012}. 
\end{proof}

For the next result, we introduce the dual problem: find $(\bzetau, \bzetap)\in \VectorH_0^1(\Omega) \times L^2(\Omega)$ such that
\begin{subequations} \label{eqs: aux Stokes}
	\begin{align}
	(\Grad \bzetau, \Grad \bfv) - (\bzetap, \Div \bfv) &= (\bfu - \bfu_h, \bfv) ,  
    &\forall \bfv&\in \VectorH_0^1(\Omega), \label{eq:Stokes1}\\
	(\Div \bzetau, q) &= 0,  
    &\forall q&\in L^2(\Omega). \label{eq:Stokes2}
	\end{align}
\end{subequations}
We are now ready to state and prove the main result of this section.

\begin{theorem} \label{thm: conv u r2}
    Let the polynomial degree $r \ge 2$. If the solution to \eqref{eqs: aux Stokes} satisfies
    \begin{align} \label{eq: bound w and phi}
        \| \bzetau \|_{2, \Omega} + \| \bzetap \|_{1, \Omega}
        \lesssim \| \bfu - \bfu_h \|_\Omega,
    \end{align} 
    then the following estimate holds
    \begin{align*}
        \lVert \bfu -\bfu_h \rVert_{\Omega} 
        &\lesssim 
        h^r \lVert \bfu \rVert_{r, \Omega} + h^{r + 1} \lVert \Curl \bfu \rVert_{r, \Omega}
        \nonumber\\
        &\quad + h \left(\inf_{\bfv_h \in \X_h}\|\bfu-\bfv_h\|_{\#} + \inf_{q_h \in \Qh}\lvert p - q_h \rvert_{1,\Omega} \right).
    \end{align*}
    In other words, the velocity converges with order $r$ in $L^2$.
\end{theorem}
\begin{proof}
By assumption, $\bzetau \in \VectorHtwo$, and so \eqref{eq:Stokes1} can be rewritten in the equivalent form 
\begin{align}
\label{eq:Stokes1a}
    (\bfv, \Curl \Curl \bzetau) + (\Grad \bzetap, \bfv) &= (\bfu - \bfu_h, \bfv),  
    &\forall \bfv&\in \VectorH_0^1(\Omega).
\end{align}
Now, since $\VectorH_0^1(\Omega)$ is dense in $\VectorLtwo$, equation \eqref{eq:Stokes1a} holds for any $\bfv \in \VectorLtwo$. In particular, taking $\bfv= \bfu - \bfu_h$ and integrating by parts, we obtain
\begin{align*}
    \lVert \bfu - \bfu_h \rVert_{\Omega}^2 
    &= (\Curl (\bfu - \bfu_h), \Curl \bzetau) + \langle \ttrace(\bfu - \bfu_h), \ntrace (\bzetau) \rangle 
    + (\Grad \bzetap, \bfu - \bfu_h) \\ 
    &= a_h(\bfu - \bfu_h, \bzetau) + b(\bfu - \bfu_h, \bzetap) \\
    &= a_h(\bfu - \bfu_h, \bzetau - \Pu) - b( \Pu, p - p_h) + b( \bfu - \bfu_h, \bzetap)
\end{align*}
in which the final equality follows from the Galerkin orthogonality property from \Cref{cor: orthogonality}. We now note that \eqref{eq:Stokes2} and \Cref{thm:Consistency} provide the identities:
\begin{align*}
    b(\bzetau, p - p_h) &= 0, &
    b( \bfu - \bfu_h, \Pp) &= 0,
\end{align*}
Using these two properties, we have
\begin{align*}
    \lVert \bfu - \bfu_h \rVert_{\Omega}^2 
    &= a_h(\bfu - \bfu_h, \bzetau - \Pu) + b(\bzetau - \Pu, p - p_h) + b(\bfu - \bfu_h, \bzetap - \Pp)
    \nonumber \\ 
    & \eqqcolon T_1 + T_2 + T_3.
\end{align*}

We bound these three terms separately.
\begin{itemize}[leftmargin=*]
    \item For $T_1$, we first use the continuity of $a_h$ from \Cref{lem: continuity}.
    \begin{align*}
        T_1 \lesssim \| \bfu - \bfu_h \|_{\#} \| \bzetau - \Pu \|_{\#}.
    \end{align*}
    We now bound the second factor by recalling the norm from \eqref{eq: hashtagnorm} and using $\ttrace(\bzetau - \Pu) = 0$.
    Note that $\lVert \ntrace(\bfw - P_{\mathring{\mathbf{V}}_h}\bfw)\rVert_{\Gamma}$ can be bound with \eqref{eq:global_trace_approx_gamman}, the inverse trace inequality \eqref{eq:trace_ineq_gamman} and \Cref{lem: projection Vzero}: 
    \begin{align*}
        \lVert \ntrace(\bfw - P_{\mathring{\mathbf{V}}_h}\bfw)\rVert_{\Gamma} &\leq \lVert \ntrace(\bfw - \mathcal{I}^c_h \bfw) \rVert_{\Gamma} + \lVert \ntrace(\mathcal{I}^c_h \bfw - P_{\mathring{\mathbf{V}}_h} \bfw)\rVert_{\Gamma} \\
        &\lesssim h^{\frac{1}{2}} \lvert \Curl \bfw \rvert_{1, \Omega} + h^{-\frac{1}{2}} \lVert\Curl( \mathcal{I}^c_h \bfw - P_{\mathring{\mathbf{v}}_h}\bfw) \rVert_{\Omega} \\
        &\leq h^{\frac{1}{2}} \lvert \Curl \bfw \rvert_{1, \Omega} + h^{-\frac{1}{2}} \lVert \Curl(\mathcal{I}^c_h \bfw - \bfw ) \rVert_{\Omega} \\
        &+ h^{-\frac{1}{2}}\lVert \Curl ( \bfw - P_{\mathring{\mathbf{V}}_h}\bfw )\rVert_{\Omega}\\
        & \lesssim h^{\frac{1}{2}} \lVert \Curl \bfw \rVert_{1, \Omega} \\
        & \leq h^{\frac{1}{2}} \lVert \bfw \rVert_{2, \Omega}.
    \end{align*}
    It follows that
    \begin{align*}
        \| \bzetau - \Pu \|_{\#} 
        &\lesssim
        \| \Curl (\bzetau - \Pu) \|_{\Omega}
        +
        h^{\frac12} \| \ntrace(\bzetau - \Pu) \|_\Gamma \\
        &\lesssim h \| \bzetau \|_{2, \Omega}
        \lesssim h \| \bfu - \bfu_h \|_\Omega
    \end{align*}
    where the final inequality is due to \eqref{eq: bound w and phi}.
    \item 
    To bound $T_2$, we again use the continuity from \Cref{lem: continuity}.
    \begin{align*}
        T_2 \leq \| \bzetau - \Pu \|_\Omega \lvert p - p_h \rvert_{1,\Omega}
    \end{align*}
    Since $r \geq 2$ by assumption, \eqref{eq:thm3.4a} from \Cref{lem: projection Vzero} gives us
    \begin{equation} \label{eq:h2}
        \lVert \bzetau - \Pu\rVert_{\Omega} 
        \lesssim h^2 \lVert \bzetau \rVert_{2, \Omega} 
        \lesssim h ^2 \lVert \bfu - \bfu_h \rVert_{\Omega}.
    \end{equation} 
    
    To conclude, we use \Cref{lem: suboptimal p conv} to bound the error in the pressure:
    \begin{align*}
        T_2 
        &\lesssim h^2\lVert \bfu - \bfu_h \rVert_{\Omega} \lvert p - p_h \rvert_{1,\Omega}
        \lesssim h\lVert \bfu - \bfu_h \rVert_{\Omega} \left( \lVert \bfu - \bfu_h \rVert_{\#} + \inf_{q_h \in \Qh}\lvert p - q_h \rvert_{1,\Omega}\right)
    \end{align*}
    \item We split $T_3$ it as 
    \begin{align} \label{eq: split T3}
        T_3 
        &= b(\bfu - P_{\Vhz} \bfu, \bzetap - \Pp) + b( P_{\Vhz} \bfu - \bfu_h, \bzetap - \Pp) \nonumber\\
        &\eqqcolon T_3' + T_3''. 
    \end{align}
    The first term can be bounded easily, using \Cref{lem: continuity,lem: projection Qh,lem: projection Vzero}, with the bound \eqref{eq: bound w and phi}:
    \begin{align*}
        T_3' 
        &\leq \lVert \Grad(\bzetap - \Pp)\rVert_{\Omega}\lVert \bfu - P_{\Vhz} \bfu \rVert_{\Omega}
        \lesssim  \lVert \bfu - \bfu_h \rVert_{\Omega}h^r \lVert \bfu \rVert_{r, \Omega}.
    \end{align*}
    For the second term, we use \eqref{eq:3.15}, \eqref{eq:thm3.4b}, and \eqref{eq: bound w and phi} to get
    \begin{align*}
    T_3'' & \lesssim h \lVert \Grad (\bzetap - \Pp) \rVert_{\Omega} \lVert \Curl (P_{\Vhz}\bfu - \bfu_h) \rVert_{\Omega} \\ 
    & \lesssim h  \lVert \bzetap \rVert_{1, \Omega} \left( \lVert \Curl (P_{\Vhz}\bfu - \bfu) \rVert_{\Omega} + \lVert \Curl (\bfu - \bfu_h) \rVert_{\Omega} \right) \\
    & \lesssim h  \lVert \bfu - \bfu_h \rVert_{\Omega} \left(h^r \lVert \Curl \bfu \rVert_{r, \Omega} + \lVert \bfu - \bfu_h \rVert_{\#}\right) .
    \end{align*}
\end{itemize}

Collecting the bounds on the three terms, we obtain
\begin{align*}
    \lVert \bfu - \bfu_h \rVert_{\Omega}
    &= \frac{T_1 + T_2 + T_3}{\lVert \bfu - \bfu_h \rVert_{\Omega}} \nonumber\\
    &\lesssim 3h \lVert \bfu - \bfu_h \rVert_{\#} + h \left(\inf_{q_h \in \Qh}\lvert p - q_h \rvert_{1,\Omega} \right)
    + h^r \lVert \bfu \rVert_{r, \Omega} + h^{r + 1} \lVert \Curl \bfu \rVert_{r, \Omega},
\end{align*}
and the result follows by applying \Cref{thm:convhashtagnorm} to the first term.
\end{proof}

\begin{theorem} \label{thm: conv u r1}
    Let the polynomial order $r = 1$ and let the solution to \eqref{eqs: aux Stokes} satisfy \eqref{eq: bound w and phi}. Set $\chi(h) \coloneqq \lvert \log h\rvert^{\frac{3}{2}}h^{\frac{1}{2}}$. Then the following estimate holds
    \begin{align*}
        \lVert \bfu -\bfu_h \rVert_{\Omega} 
        &\lesssim 
        h(1+\chi(h)) \lVert p \rVert_{1, \Omega}
        + h \lVert \bfu \rVert_{1, \Omega} + h^2 \lVert \Curl \bfu \rVert_{1, \Omega} \nonumber\\
        &\quad + (2h + \chi(h))\left(\inf_{\bfv_h \in \X_h}\|\bfu-\bfv_h\|_{\#} + \inf_{q_h \in \Qh}\lvert p - q_h \rvert_{1,\Omega}\right).
    \end{align*}
    
    Since $\chi(h)h^{1/2}\lesssim O(h)$ when $h$ goes to zero, in the asymptotic limit the velocity converges linearly in $L^2$. 
\end{theorem}
\begin{proof}
    Following the proof of \Cref{thm: conv u r2}, we only used $r \ge 2$ in the bound for $T_2$. In particular, estimate \eqref{eq:h2} does not hold for $r = 1$, so we require a different approach. Since the preceding steps are the same, we directly consider the bound on $T_2$. 
    \begin{itemize}[leftmargin=*]
        \item 
        We follow the proof of \cite[Thm.~3.9]{ArFaGo2012}. Similar to \eqref{eq: split T3}, we split $T_2$ as
        \begin{align} \label{eq: split T2}
            T_2 &= b(  \bzetau - \Pu, p - P_{\Qh}p) + b(  \bzetau - \Pu, P_{\Qh}p - p_h) \nonumber\\
            &\eqqcolon T_2' + T_2''.
        \end{align} 
        
        Bounding the first term is straightforward. In particular, continuity of $b$, the approximation properties \eqref{eq:thm3.4a} (with $r = 1$) and \Cref{lem: projection Qh}, and the bound \eqref{eq: bound w and phi} give us
        \begin{align*}
            T_2' \leq \lVert \bzetau - \Pu\rVert_{\Omega} \lVert \Grad ( p - P_{\Qh}p)\rVert_{\Omega} 
            &\lesssim h\lVert \bzetau \rVert_{1, \Omega} \lVert p\rVert_{1, \Omega} \\
            &\lesssim h \lVert \bfu - \bfu_h \rVert_{\Omega} \lVert p\rVert_{1, \Omega}
        \end{align*}

        The second term, on the other hand, is bound using \eqref{eq:thm3.5} and \eqref{eq:thm3.4a} from \Cref{lem: projection Vzero}:
        \begin{align*}
            T_2''
            & \lesssim h^{-\frac{1}{2} - \frac{1}{q}}\lVert \bzetau - \Pu\rVert_{\Omega,q} \lVert P_{\Qh}p - p_h \rVert_{\Omega} \\
            & \lesssim qh^{\frac{1}{2} - \frac{1}{q}}\lVert \bzetau\rVert_{1,\Omega, q} \lVert P_{\Qh} p - p_h\rVert_{\Omega},
        \end{align*}
        for each $ 2\leq q < \infty$. Since 
        $\lVert \bzetau \rVert_{1,\Omega,q} \leq  Cq^{\frac{1}{2}} \lVert \mathbf{w} \rVert_{2, \Omega}$ \cite[Thm.~3.7]{ArFaGo2012}
        , we have
        \begin{equation*}
            T_2'' \lesssim q^{\frac{3}{2}}h^{\frac{1}{2} - \frac{1}{q}} \lVert \bfu - \bfu_h \rVert_{\Omega} \lVert P_{\Qh} p - p_h \rVert_{\Omega}
        \end{equation*}
        It remains to estimate the pressure term, which we do by employing \Cref{lem: projection Qh,lem: p conv L2}
        \begin{align*}
            \lVert P_{\Qh} p - p_h \rVert_{\Omega} 
            & \leq \lVert p - P_{\Qh} p \rVert_{\Omega} + \lVert p - p_h \rVert_{\Omega} \\
            & \lesssim h\lVert p \rVert_{1, \Omega} + \lVert \bfu - \bfu_h \rVert_{\#} + \inf_{q_h \in \Qh}\lvert p - q_h \rvert_{1,\Omega}.
        \end{align*}
    \end{itemize}
    Taking $q = \lvert \log h\rvert$, using $h^{ - \frac{1}{\lvert \log h \rvert}}=e$ if $h<1$, and combining the bounds for $T_2'$ and $T_2''$, we obtain 
    \begin{align*}
        T_2 \lesssim h \lVert \bfu -\bfu_h \rVert_{\Omega} + h\chi(h) \lVert p \rVert_{1} + \chi(h)\left(\lVert \bfu - \bfu_h \rVert_{\#} + \inf_{q_h\in \Qh}\lVert  p - q_h \rVert_{Q} \right).
    \end{align*}
    The rest of proof is identical to the one of \cref{thm: conv u r2}.
\end{proof}
\begin{remark}\label{rmk:L2_3D}
In \Cref{thm: conv u r2} and \Cref{thm: conv u r1}, we assume that $\Omega$ is two-dimensional, as the argument depends on the existence of the projection operator $P_{\mathring{\mathbf{V}}_h}$. If, however, a projection operator with analogous properties exists in three dimensions, the results of \Cref{thm: conv u r2} and \Cref{thm: conv u r1} can be extended to the three-dimensional setting in a straightforward manner.
\end{remark}

\section{Numerical experiments}
\label{sec:numerics}
In this section, we present a series of numerical experiments to validate the theoretical results established in the preceding chapters. In addition, we examine whether the assumptions in \Cref{thm: conv hashtag norm}, \Cref{thm: conv u r2}, and \Cref{thm: conv u r1} are strictly necessary. As anticipated in \Cref{rem:improvedconvergence3D} and \Cref{rmk:L2_3D}, the numerical results support the conjecture that the improved convergence rates extend to three-dimensional and non-topologically trivial domains.  For clarity and comparison, the theoretical and observed convergence rates are summarized in \Cref{tab:convergence_rates}. All simulations were carried out using the NGSolve library \cite{ngSolve}. The source code is freely available at \url{https://gitlab.com/WouterTonnon/vvphcurlslip}.

\begin{table}[ht]
\centering
\caption{Predicted and observed convergence rates in \cref{sec:manufacturedSolutionStarShaped2D,sec:domainwithhole,sec:ManufacturedSolution3D,sec:FlowAroundSlipperySphere} }
\small
\begin{tabular}{|l|c|cc|cc|}
\cline{3-6}
 \multicolumn{2}{c|}{ }& 
 \multicolumn{2}{|c|}
 {\textbf{\Cref{sec:manufacturedSolutionStarShaped2D}}}& \multicolumn{2}{|c|}{\textbf{\Cref{sec:domainwithhole,sec:ManufacturedSolution3D,sec:FlowAroundSlipperySphere}}} 
 \\ \hline
 \textbf{Norm} & \textbf{Theorem} 
 & Pred. & Obs. & Pred. & Obs. \\ \hline
    $\lVert \bfu-\bfu_h\rVert_\Omega$ &\ref{thm: conv hashtag norm}, \ref{thm: conv u r2}, \ref{thm: conv u r1}&
    $r$ &$r$&  N/A & $r$\\ 
    $\lVert \Curl\bfu-\Curl\bfu_h\rVert_\Omega$ &\ref{thm: conv hashtag norm}&
    $r-\frac{1}{2}$ &$r-\frac{1}{2}$& $r-1$ & $r-\frac{1}{2}$ \\ 
    $\lVert p-p_h\rVert_\Omega$ &\ref{thm: conv hashtag norm}, \ref{lem: p conv L2}&
    $r-\frac{1}{2}$ &$r-\frac{1}{2}$& $r-1$ & $r-\frac{1}{2}$ \\ 
    $\lVert \nabla p - \nabla p_h\rVert_\Omega$ &\ref{thm: conv hashtag norm} &
    $r-\frac{3}{2}$ &$r-\frac{3}{2}$& $r-2$ & $r-\frac{3}{2}$ \\ \hline
\end{tabular}
\label{tab:convergence_rates}
\end{table}

\subsection{Convergence analysis on a star-shaped domain} \label{sec:manufacturedSolutionStarShaped2D}
\begin{figure}[htb]
  \centering
  \includegraphics[width = 0.49\textwidth]{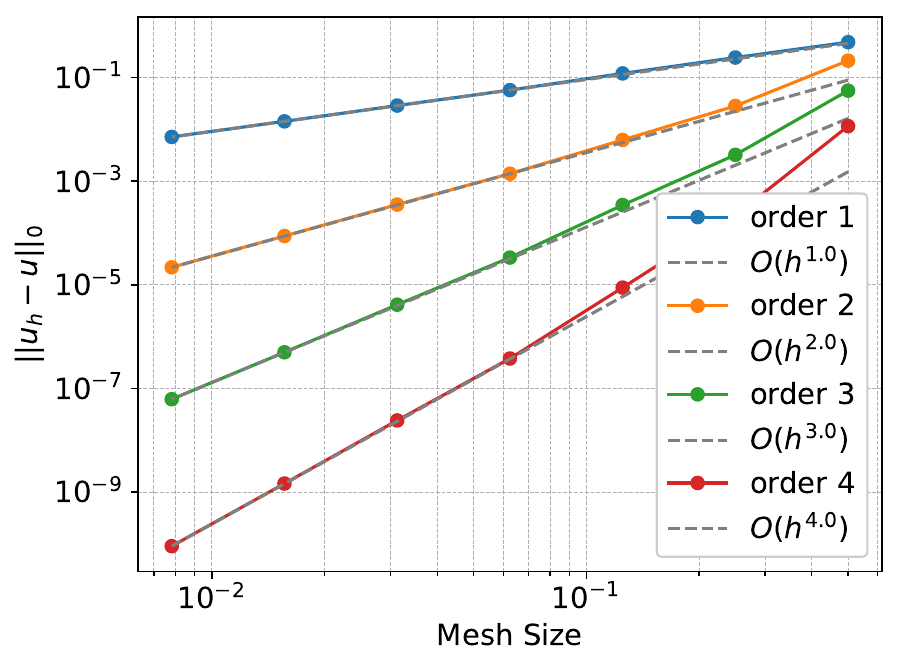}
  \includegraphics[width = 0.49\textwidth]{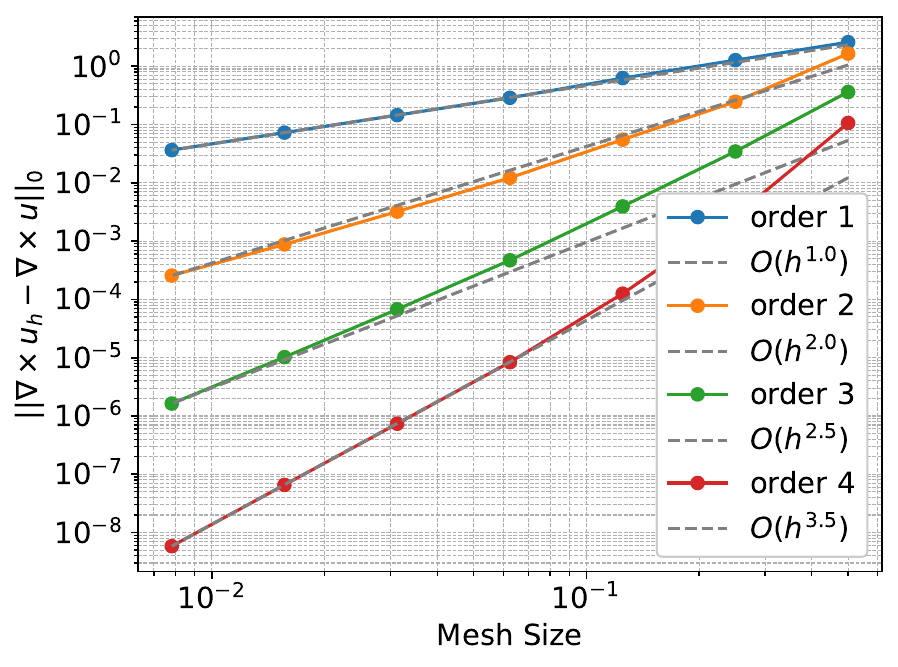} \\
  \includegraphics[width = 0.49\textwidth]{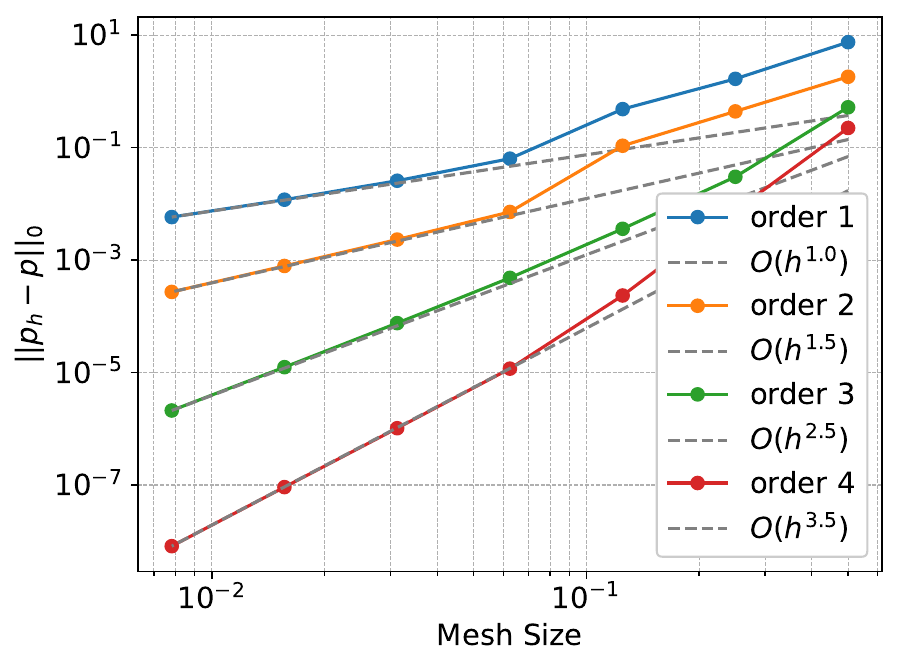}
  \includegraphics[width = 0.49\textwidth]{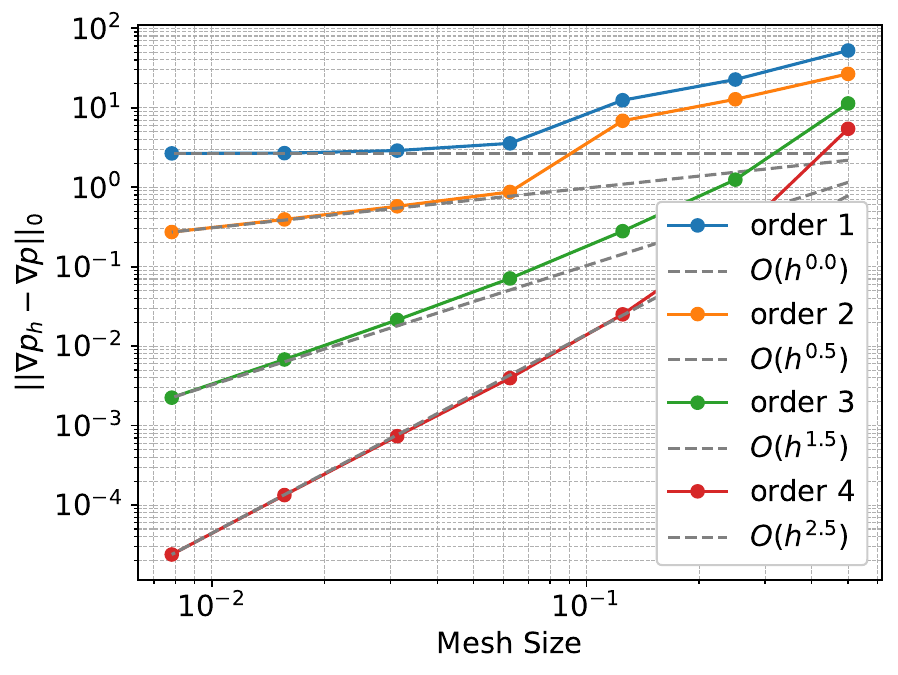}
  \caption{Convergence analysis of (top) $u$ in the $L^2$ and $\Hcurl$ norms, (middle) $\bfu$ and $\curl\bfu$ in the $L^2$ norm on the boundary, and (bottom) $p$ in the $L^2$ and $H^1$ norms for the experiment as discussed in \cref{sec:manufacturedSolutionStarShaped2D}. The results for lowest-order elements are labeled as \lq\lq order 1\rq\rq. 
  }
  \label{fig:manufacturedSolutionStarShaped2D}
\end{figure}

In this experiment, we compare the theoretical and experimental convergence on a domain that does not allow for harmonic forms. We approximate the solution to \cref{eq:continuousSymGradStokes} on a unit square with unstructured meshes. The force term $\bff$ and the boundary term $\bfg$ have been chosen in such a way that the solution is
\begin{align*}
    \bfu &= \begin{bmatrix}-\sin(4x)\cos(4y)\\ \cos(4x)\sin(4y)\end{bmatrix},\qquad
    p = \cos(4\pi x)+ \cos(4\pi y).
\end{align*}
In \cref{fig:manufacturedSolutionStarShaped2D} we display various errors and how they relate to the size of the mesh. We summarize the observed convergence rates in \cref{tab:convergence_rates}. We observe that the experimental convergence agree with the theory presented in this work.
\subsection{Convergence analysis on a domain with a hole} \label{sec:domainwithhole}
\begin{figure}[htb]
  \centering
  \includegraphics[width = 0.49\textwidth]{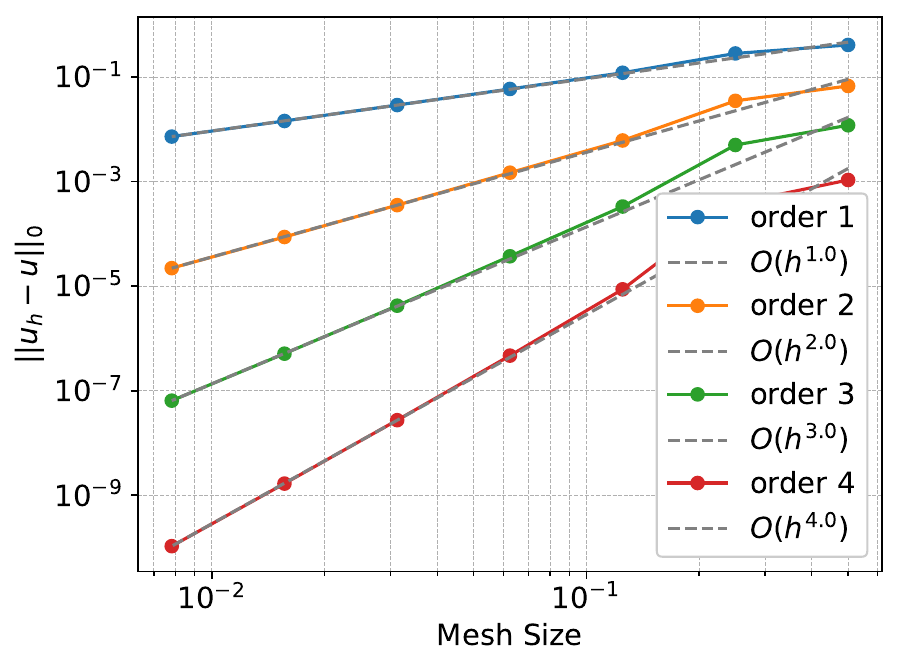}
  \includegraphics[width = 0.49\textwidth]{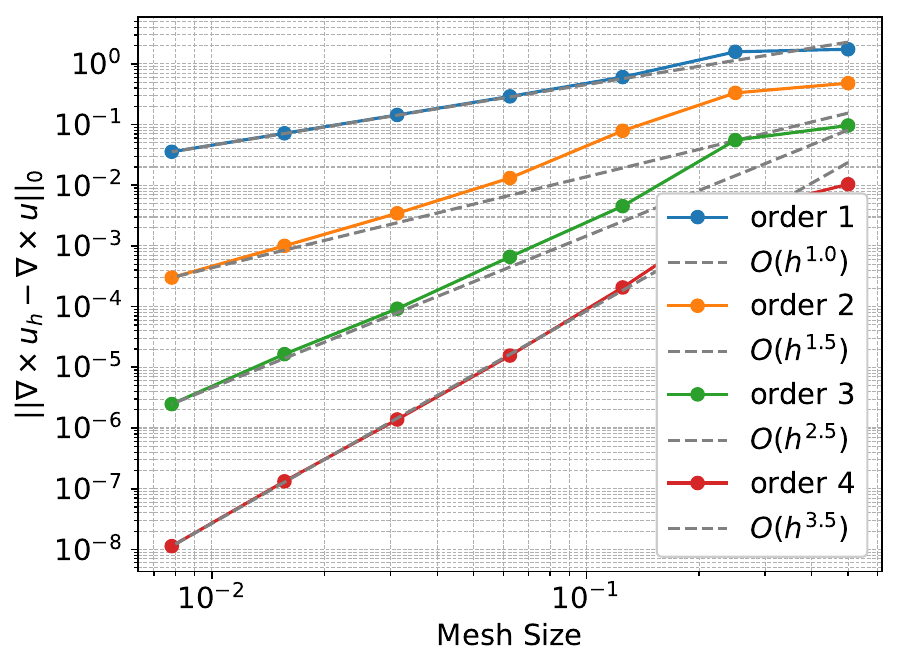} \\
  \includegraphics[width = 0.49\textwidth]{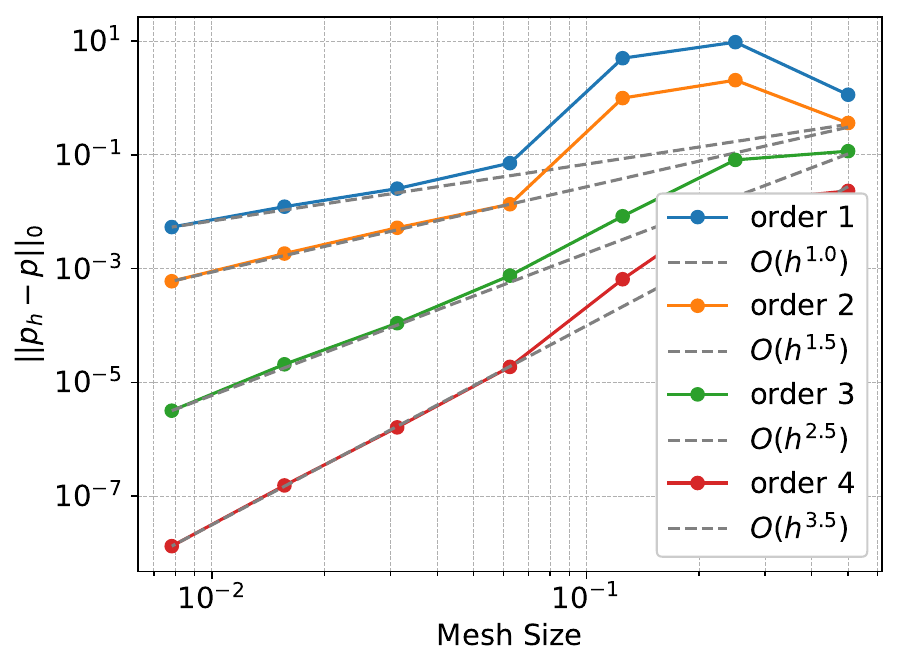}
  \includegraphics[width = 0.49\textwidth]{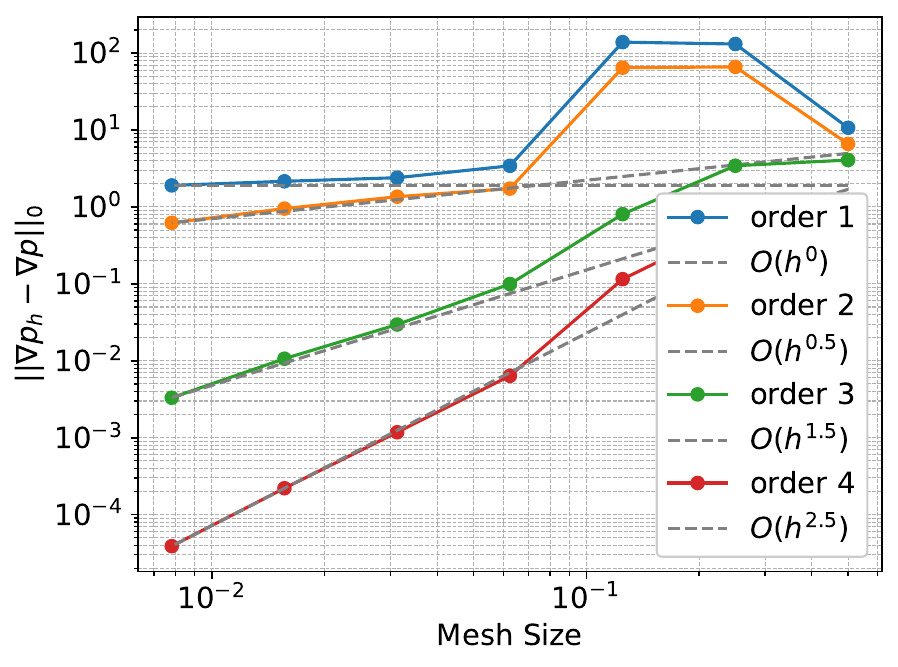}
  \caption{Convergence analysis of (top) $u$ in the $L^2$ and $\Hcurl$ norms, and (bottom) $p$ in the $L^2$ and $H^1$ norms for the experiment as discussed in \cref{sec:domainwithhole}. The results for lowest-order elements are labeled as \lq\lq order 1\rq\rq. 
  }
  \label{fig:convergencewithhole}
\end{figure}

We repeat the previous experiment, but this time the domain is $[0,1]^2 \setminus [\frac{1}{3},\frac{2}{3}]^2$, which allows for harmonic forms and has reentrant corners. We approximate the solution to \cref{eq:continuousSymGradStokes} on . 
In \cref{fig:convergencewithhole} we display various errors and how they relate to the size of the mesh. 
\subsection{A manufactured solution in 3D} \label{sec:ManufacturedSolution3D}
\begin{figure}[htb]
  \centering
  \includegraphics[width = 0.49\textwidth]{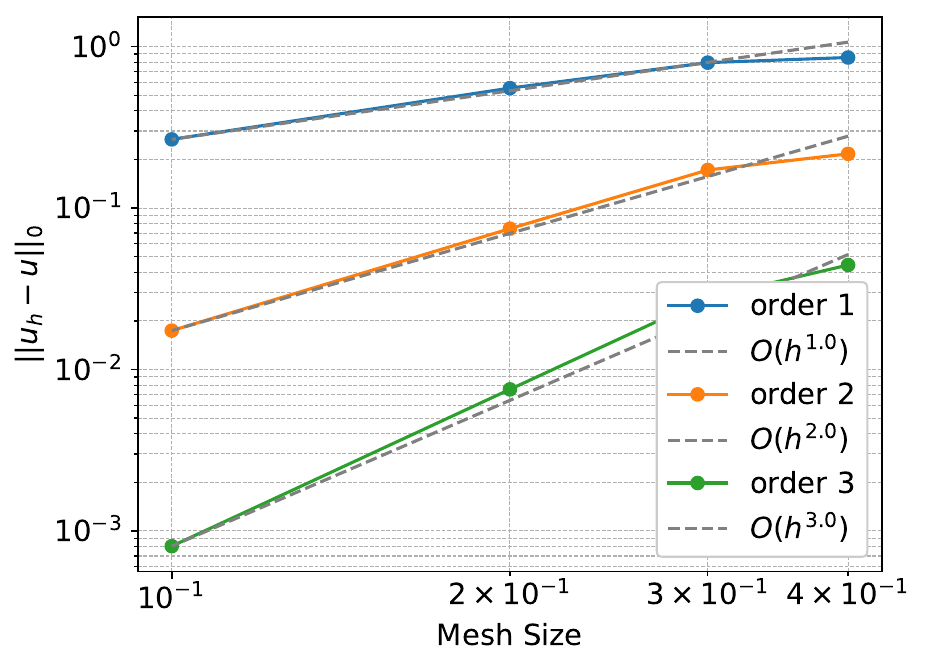}
  \includegraphics[width = 0.49\textwidth]{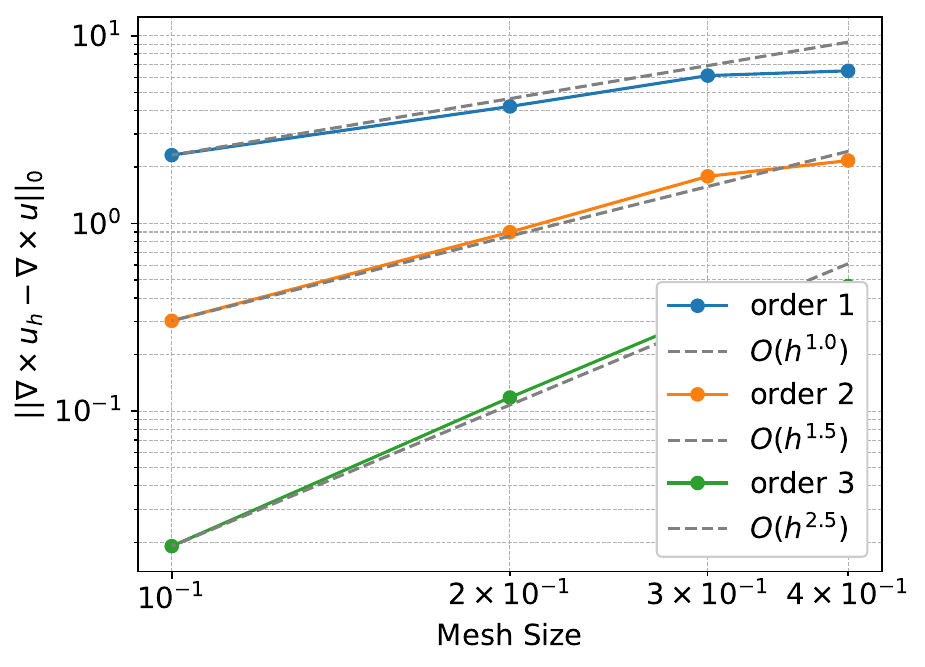} \\
  \includegraphics[width = 0.49\textwidth]{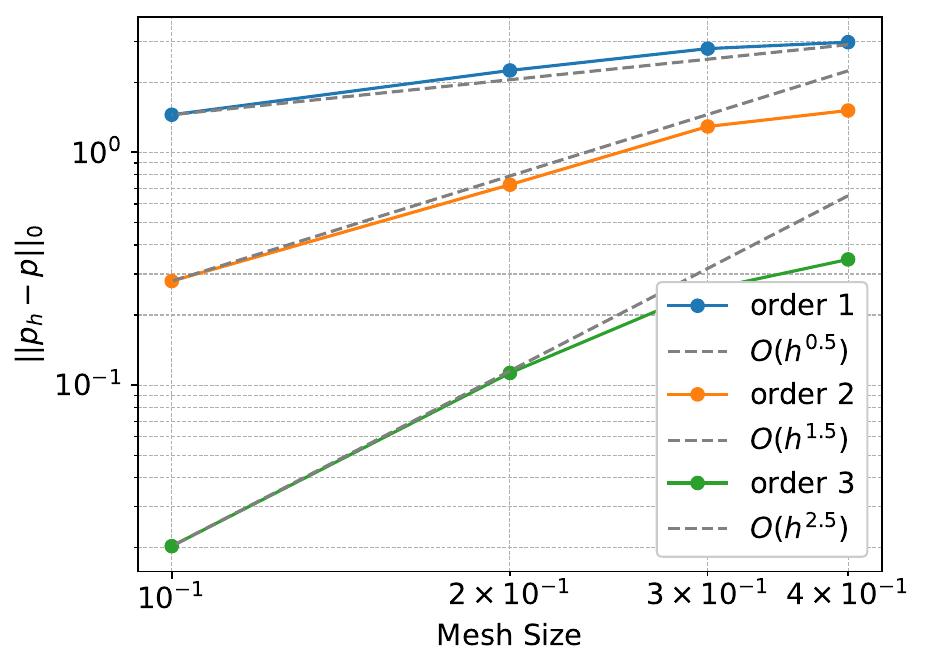}
  \includegraphics[width = 0.49\textwidth]{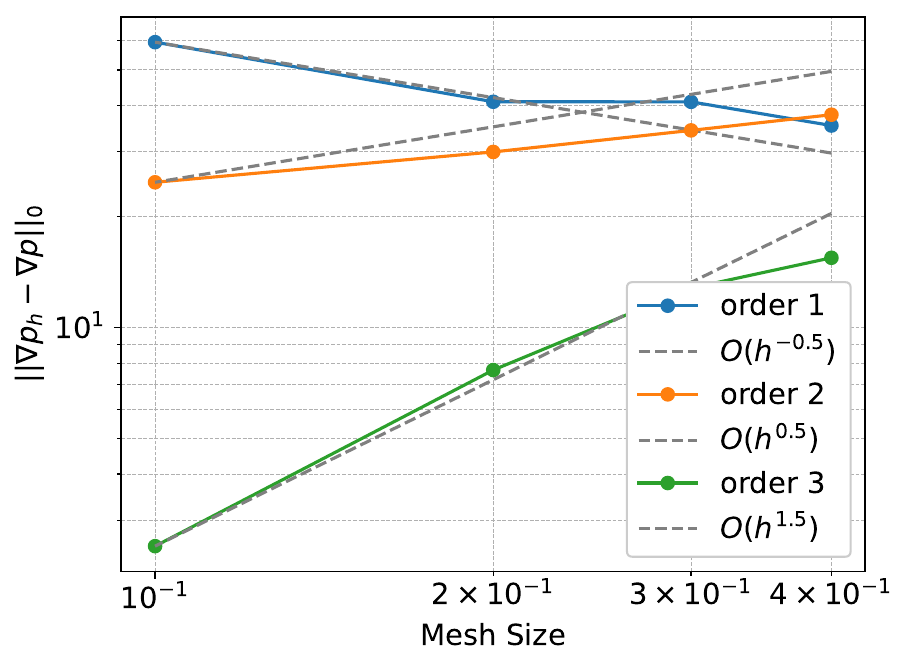}
  \caption{{Convergence analysis of (top) $u$ in the $L^2$ and $\Hcurl$ norms, and (bottom) $p$ in the $L^2$ and $H^1$ norms for the experiment as discussed in \cref{sec:ManufacturedSolution3D}. The results for lowest-order elements are labeled as \lq\lq order 1\rq\rq. 
  }}
  \label{fig:NitscheConvergenceAnalysis3D}
\end{figure}

In this experiment, we consider a manufactured solution on an ellipsoid $\Omega\subset\R^3$ with width 1, height 0.8, and depth 1.2. Given the solution
\begin{align*}
    \bfu &= \Curl\begin{bmatrix}
       \sin(\pi y)\cos(\pi z)x^2 \\
       \sin(\pi x)\cos(\pi x)yz\\
       \sin(\pi z)\cos(\pi x)xz^2
    \end{bmatrix}, &
    p = \frac{1}{10}\sin(\pi x)\cos(\pi y)\cos(\pi z),
\end{align*}
we derive $\mathbf{f}$, $z$, and $\mathbf{g}$. In \Cref{fig:NitscheConvergenceAnalysis3D}, we display the $L^2$-error of $\bfu$, the $\Hcurl$-error of $\bfu$, the $L^2$-error of $p$, and the $H^1$-error of $p$, respectively. 
\subsection{Flow around sphere}\label{sec:FlowAroundSlipperySphere}
{
In this experiment, we consider Stokes flow around a sphere in 3D, following \cite{CostaSlip}. The domain is a 2-by-2-by-2 (bounding) box with a sphere of radius 0.5 at the center cut out from the box. The exact solution of this problem in polar coordinates is 
\begin{equation*}
    u_r = \cos(\theta)\left(1+\frac{a^3}{2r^3}-\frac{3a}{4r}\right), \quad u_{\theta} = \sin(\theta)\left(1-\frac{a^3}{4r^3} - \frac{3a}{4r} \right), \quad p = \frac{3a}{2r^3}\cos(\theta),
\end{equation*}
where $r$ indicates the radius and $\theta$ indicates the polar angle. Note that the solution is axisymmetric around the $z$-axis and, thus, is independent of the azimuthal angle.
We visualize the results in \cref{fig:NitscheFlowCylinderConvergenceAnalysis3D}.
We also perform a convergence analysis and report the result in \Cref{fig:NitscheFlowCylinderConvergenceAnalysis3D}.
\begin{figure}[htb]
  \centering
  \includegraphics[width = 0.3\textwidth]{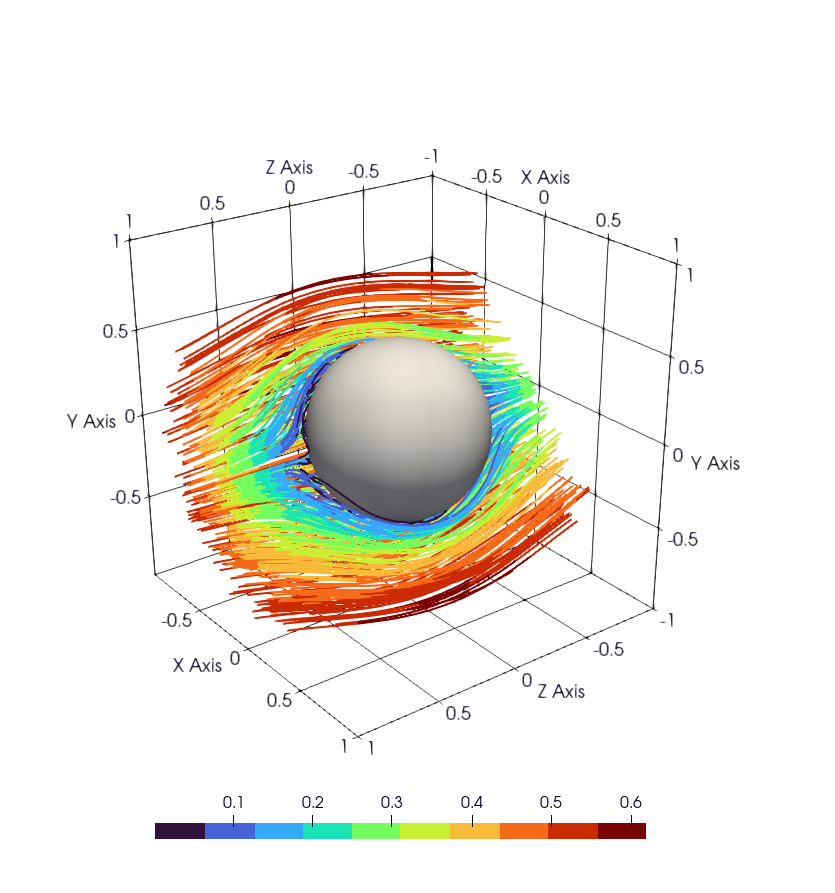}
  \includegraphics[width = 0.3\textwidth]{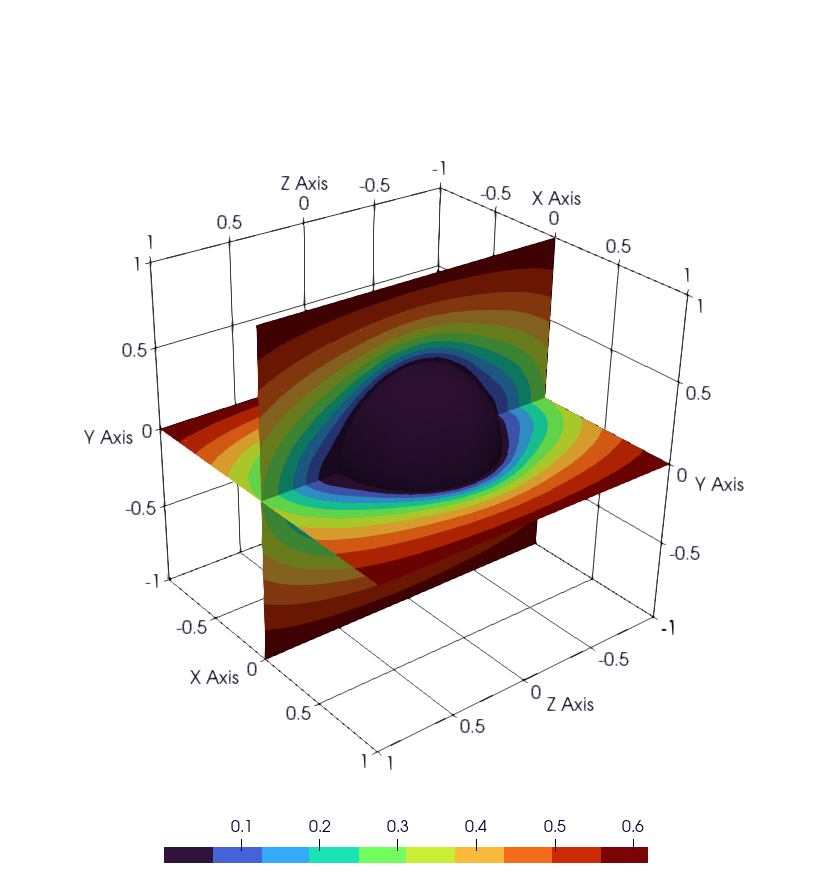}
  \includegraphics[width = 0.3\textwidth]{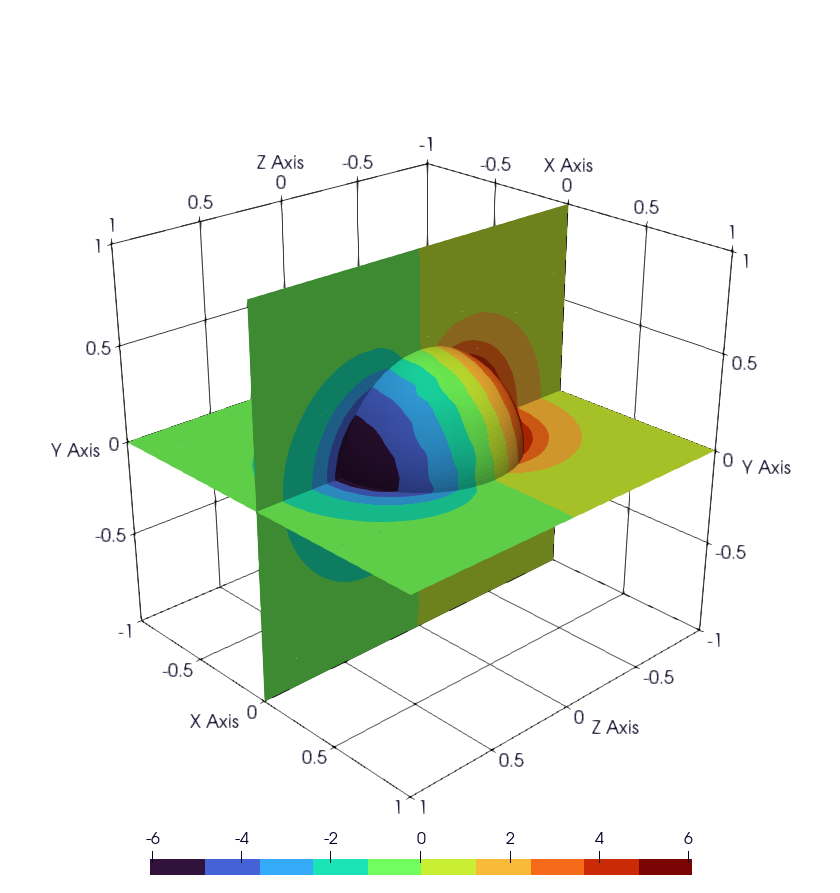}
  \caption{{Visualization of the computed solution as described in \cref{sec:FlowAroundSlipperySphere}. The domain is a 2-by-2-by-2 box with a sphere of radius 0.5 cut out at its center. We used 3rd order polynomials on an unstructured, curved mesh with mesh-width $h=0.2$. (Left) The lines represent the streamlines and the colors indicate the magnitude of the velocity field $\bfu$. (Middle) The colors indicate the magnitude of the velocity field. (Right) The colors indicate the pressure.}}
  \label{fig:FlowAroundSphere}
\end{figure}
\begin{figure}[ht]
  \centering
  \includegraphics[width = 0.49\textwidth]{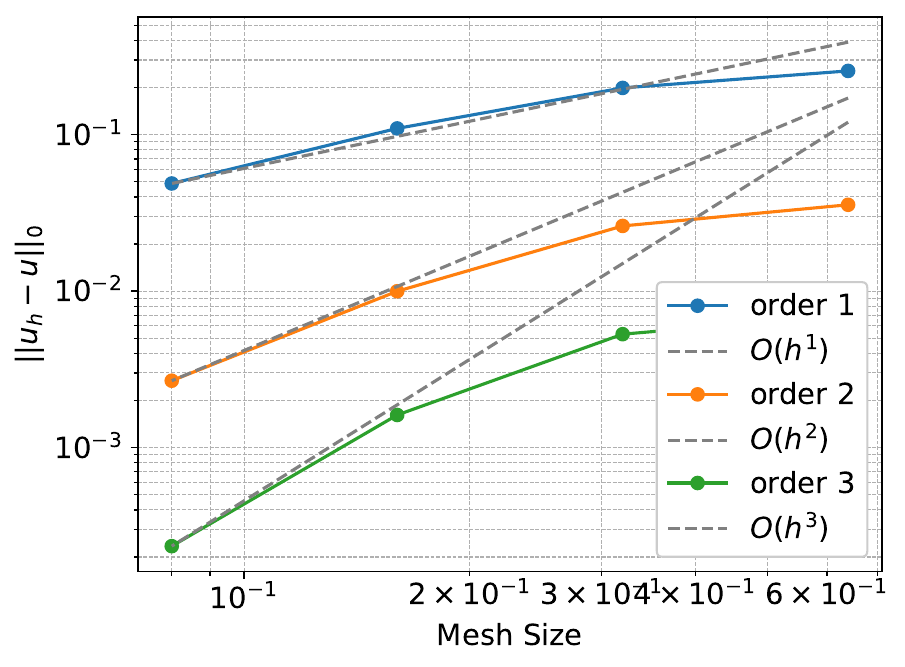}
  \includegraphics[width = 0.49\textwidth]{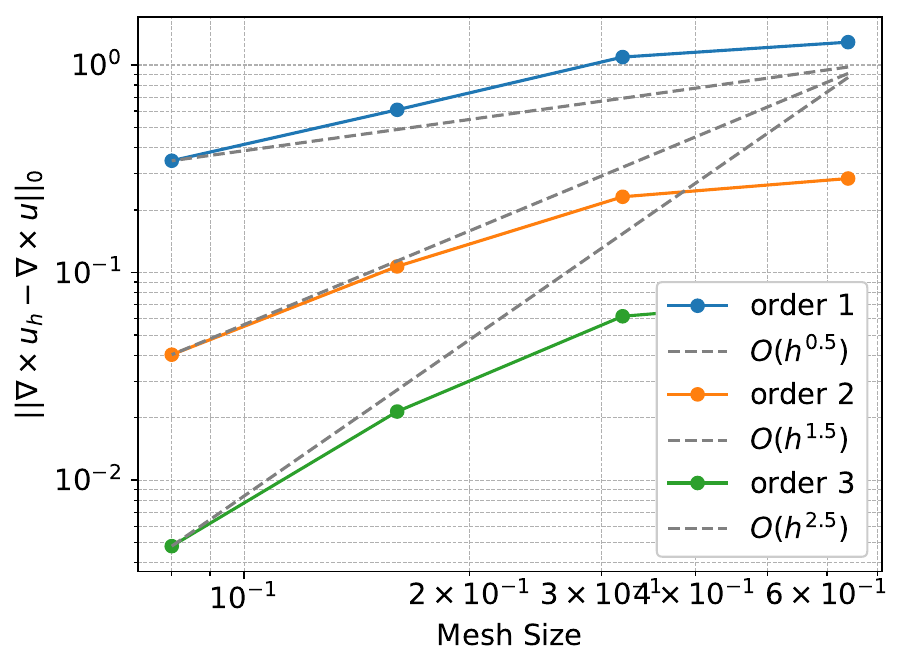} \\
  \includegraphics[width = 0.49\textwidth]{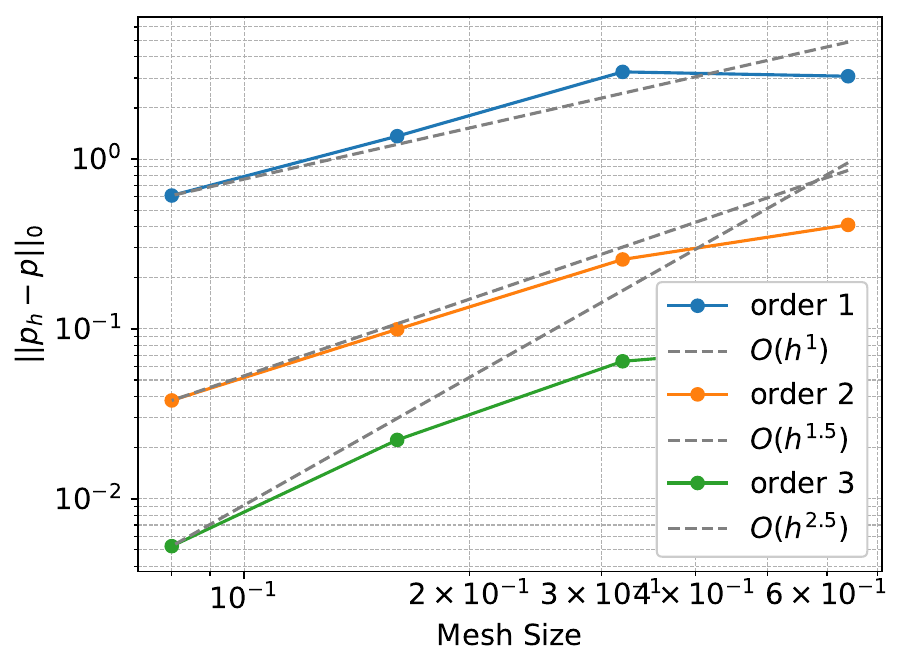}
  \includegraphics[width = 0.49\textwidth]{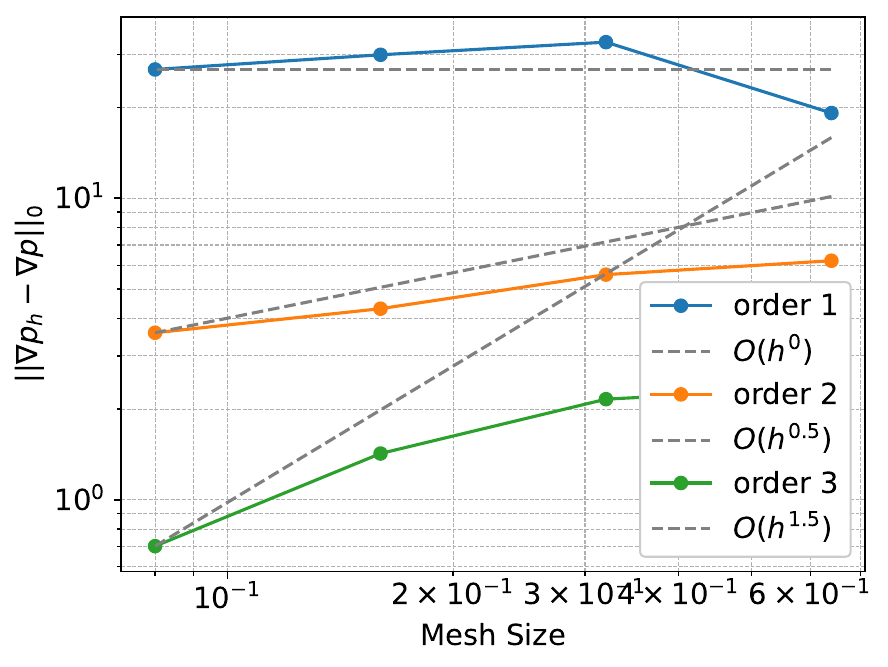}
  \caption{{Convergence analysis of (top) $u$ in the $L^2$ and $\Hcurl$ norms, and (bottom) $p$ in the $L^2$ and $H^1$ norms for the experiment as discussed in \cref{sec:FlowAroundSlipperySphere}. The results for lowest-order elements are labeled as \lq\lq order 1\rq\rq. 
  }}
\label{fig:NitscheFlowCylinderConvergenceAnalysis3D}
\end{figure}

\subsection{A non-smooth solution} \label{sec:NonSmoothExperiment}

To investigate if the smoothness assumptions in \Cref{thm:Consistency} are strictly necessary,  we consider a solution on an L-shaped domain $\Omega= (-1,1)^2 \setminus [0,1)\times (-1,0]$  with $\bff=\mathbf{0}$, following \cite{HoustonSchoetzau,verfuhrt1996review}. Let $(r,\phi)$ denote the standard polar coordinates, then we seek the solution
\begin{align*}
    \bfu &= \begin{bmatrix}
        r^\lambda(1+\lambda)\sin(\phi)\Psi(\phi)+\cos(\phi)\Psi''(\phi)\\
        r^\lambda\sin(\phi)\Psi'(\phi)-(1+\lambda)\cos(\phi)\Psi(\phi)
    \end{bmatrix}, \\
    p &= -r^{\lambda-1}[(1+\lambda)^2\Psi'(\phi)+\Psi'''(\phi)]/(1-\lambda),
\end{align*}
where
\begin{align*}
    \Psi(\phi) &= \sin((1+\lambda)\phi)\cos(\lambda\omega)/(1+\lambda)-\cos((1+\lambda)\phi)\\
    &\quad -\sin((1-\lambda)\phi)\cos(\lambda\omega)/(1-\lambda)+\cos((1-\lambda)\phi),\\
    \omega &=\frac{3\pi}{2},
\end{align*}
and $\lambda\approx 0.54448373678246$. Note that both $\nabla\bfu$ and $\nabla p$ are singular at the origin, in particular $\bfu\not\in\VectorHtwo$ and $p\notin H^1(\Omega)$.  

We display the $L^2$-error of $\bfu$, the $\Hcurl$-error of $\bfu$, the $L^2$-error of $p$ in \Cref{fig:NitscheConvergenceAnalysisLshape}. We do observe convergence in all norms but at a limited rate due to the lack of regularity of the solution.
\begin{figure}[htb]
  \centering
  \includegraphics[width = 0.49\textwidth]{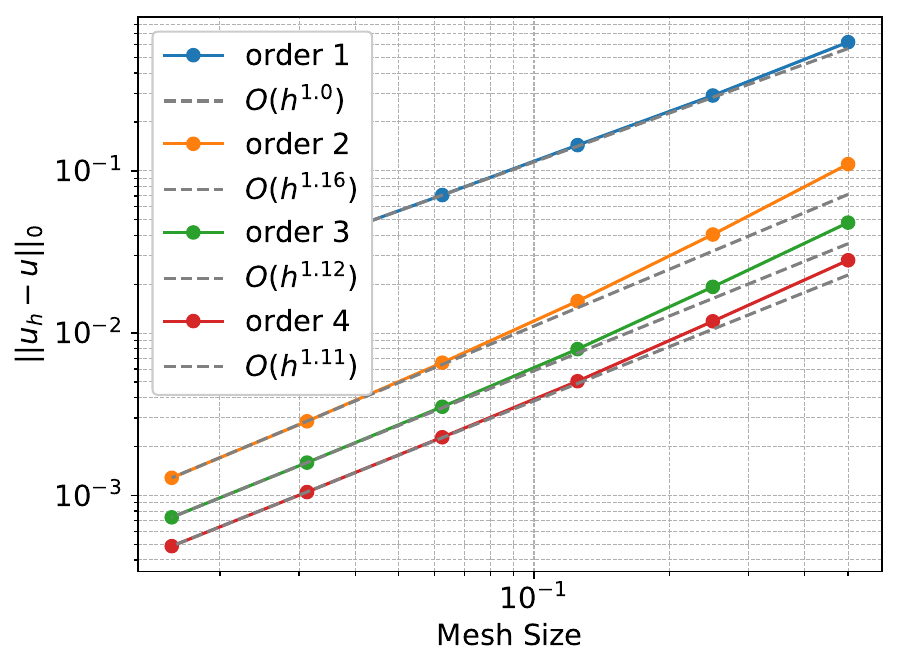}
  \includegraphics[width = 0.49\textwidth]{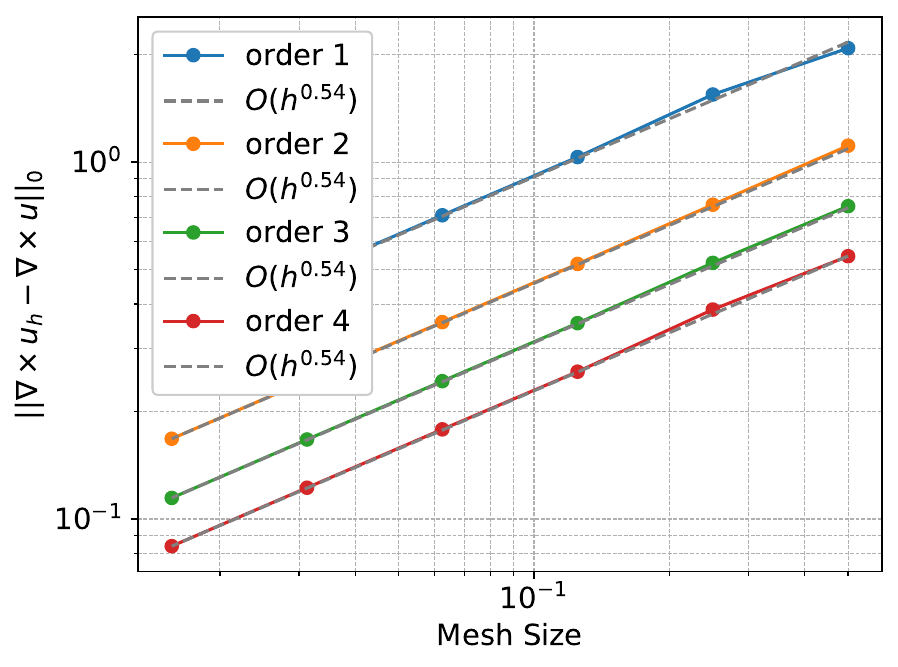} \\
  \includegraphics[width = 0.49\textwidth]{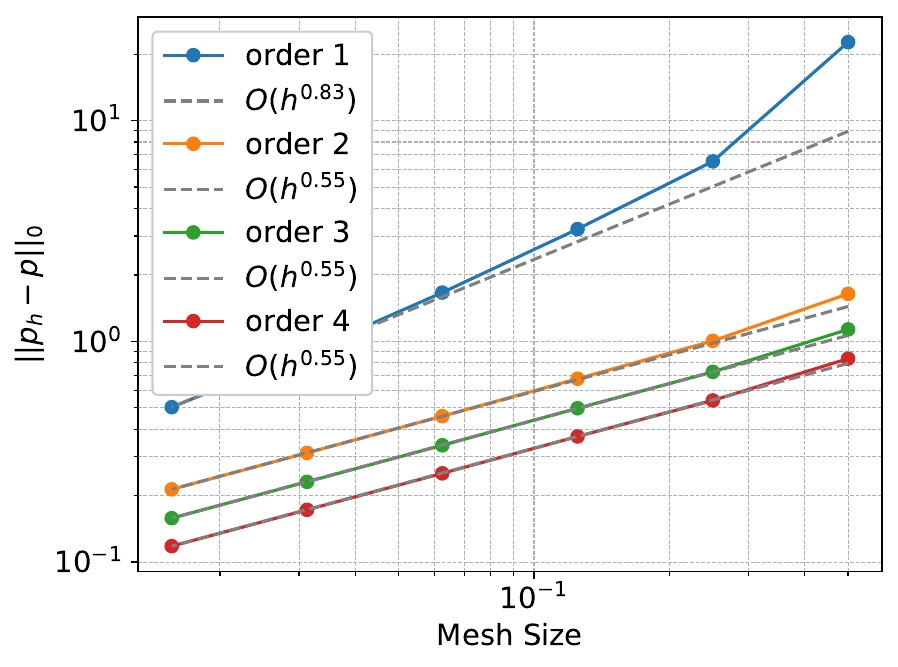}
  \caption{{Convergence analysis of (top) $u$ in the $L^2$ and $\Hcurl$ norms, and (bottom) $p$ in the $L^2$ and $H^1$ norms for the experiment as discussed in \cref{sec:NonSmoothExperiment}. The results for lowest-order elements are labeled as \lq\lq order 1\rq\rq.}}
  \label{fig:NitscheConvergenceAnalysisLshape}
\end{figure}

\section{Concluding remarks}
\label{sec:concluding_remarks}

In the context of magnetohydrodynamical (MHD) systems, it is convenient to seek the fluid velocity in $\Hcurl$, so that cross-helicity can be preserved. However, we have shown that imposing no-slip conditions as essential conditions on this space leads to an ill-posed problem. To remedy this, we have formulated and analyzed a Nitsche-type approach for the weak imposition of no-slip boundary conditions, on the simplified case of Stokes flow. 

The additional terms introduced by Nitsche's method are not continuous in the standard $\Hcurl$ norm and we therefore introduced mesh-dependent norm for the velocity. This resulted in an inf-sup constant that depends on the mesh size $h$, which is reflected in the stability estimate for the pressure. Consequently, our a priori error estimates predicted a convergence loss of at least half an order in both variables in the mesh-dependent norms. We then improved these estimates in $L^2$ using duality techniques. The predicted stability and convergence of the method was confirmed by four numerical experiments. 

In summary, we proposed a viable approach to impose no-slip boundary conditions on $\Hcurl$-based approximations of Stokes-type problems. The method is backed by rigorous analysis, is easily implementable, and results in only a minor loss in convergence for the pressure variable. 

\section*{Acknowledgments}
The authors would like to thank Ralf Hiptmair for his constructive feedback and helpful suggestions. E.Z. is member of the GNCS-INdAM (Istituto Nazionale di Alta Matematica) group.

\bibliographystyle{plain}
\bibliography{references}

\end{document}